\documentclass{article}
\usepackage{amssymb}
\newcommand{\eq}{\begin{equation}}
\newcommand{\en}{\end{equation}}
\newcommand{\te}{\rightarrow}
\newcommand{\re}[1]{\mbox{(\ref{#1})}}

\newcommand{\prob}{\mathbb P}
\newcommand{\ulambda}{\underline{n}}
\newcommand{\nlambda}{n}
\renewcommand{\ulambda}{n_1, \ldots, n_k}
\newcommand{\NLambda}{N}
\newcommand{\kell}{k}
\newcommand{\ex}{\mathbb E}
\newcommand{\nints}{\mathbb N}

\newcommand{\ed}{ \stackrel{d}{=}}
\newcommand{\tOm}{\widetilde{\Om}} 
\newcommand{\Om}{\nu}
\newcommand{\OW}{\Phi}
\newcommand{\ow}{\Phi}
\newcommand{\drift}{{\tt d}}
\def\endpf{\hfill $\Box$ \vskip0.5cm}
\def \proof{\noindent{\it Proof.\ }}

\newtheorem{theorem}{\large Theorem}[section]
\newtheorem{proposition}[theorem]   {\large Proposition}
\newtheorem{definition}[theorem]{\large Definition}
\newtheorem{corollary}[theorem]{\large  Corollary}
\newtheorem{lemma}[theorem]{\large  Lemma}

\begin{document}

\title{Regenerative Composition Structures
\thanks{Research supported in part by N.S.F. Grant DMS-0071448}
}
\author{Alexander Gnedin\thanks{Utrecht University; e-mail gnedin@math.uu.nl}
\hspace{.2cm}
and 
\hspace{.2cm}
Jim  Pitman\thanks{University of California, Berkeley; e-mail pitman@stat.Berkeley.EDU}  
\\
\\
\\
\\
}
\date{}
\maketitle

\centerline{\bf Abstract}

A new class of random composition structures (the ordered analog of Kingman's 
partition structures) is defined by a 
regenerative description of component sizes. Each regenerative composition 
structure is represented by a process of random sampling of points from an exponential distribution on the positive halfline, 
and separating the points into clusters by an independent regenerative 
random set.
Examples are composition structures derived from residual allocation models, including one associated with the Ewens sampling formula,
and composition structures derived from the zero set of a Brownian motion or
Bessel process.
We provide characterisation results and formulas relating the distribution of the regenerative
composition to the L{\'e}vy parameters of a subordinator whose range is the corresponding regenerative set.
In particular, the only reversible regenerative composition 
structures are those associated with the interval partition of $[0,1]$ generated by excursions of a standard Bessel bridge of dimension 
$2 - 2 \alpha$ for some $\alpha \in [0,1]$.

\vskip0.5cm

{\it AMS} 2000 {\it subject
classifications.} Primary 60G09, 60C05. Keywords: 
exchangeability,
composition structure, regenerative set, sampling formula, subordinator

\vskip0.5cm

\newpage

\section{Introduction}

A {\it composition} of a positive integer $n$ is a
sequence of positive integers $\lambda=(\nlambda_1,\ldots,\nlambda_{\kell})$ with
sum $\Sigma_j n_j=n$.
Each $\nlambda_i$ may be called a 
{\em part} of the composition. 
We will use the notation $\lambda\models n$ to say that $\lambda$ is a composition of $n$.
A {\em random composition} of $n$ is a random variable
${\cal C}_n$ with values in the set 
of all $2^{n-1}$ compositions of $n$.
A {\it composition structure} $({\cal C}_n)$ is a Markovian sequence of 
random compositions of $n$, for $n = 1,2, \ldots$,
whose cotransition probabilities are determined by the
following 
property of {\it sampling consistency} \cite{dj91, gnedin97}:
if $n$ identical balls are distributed
into an ordered series of boxes according to 
$({\cal C}_n)$, then
${\cal C}_{n-1}$ is obtained by discarding
one of the balls picked uniformly at random, and then deleting an 
empty box in case one is created.  We study composition structures with the following further property:

\begin{definition}
\label{regendef}
{\em A composition structure $({\cal C}_{n})$ is {\em regenerative}
if for all $n>m\geq 1$, 
given that the first part of ${\cal C}_{n}$ is $m$,
the remaining composition of $n-m$ is distributed like ${\cal C}_{n-m}$.
}
\end{definition}
According to our main result (Theorem \ref{thm1}),
each regenerative composition structure
can be represented by a process of random sampling of points from the exponential distribution on $[0,\infty[$,
and separating the sample points into clusters by points of an independent regenerative random closed subset  $\cal R$ of $[0,\infty[$.
We recall in Theorem \ref{maison} the fundamental result of 
Maisonneuve \cite{mais83} that every such $\cal R$ can
be represented as the closed range of a {\em subordinator} $(S_t)$, that
is an increasing process with stationary independent increments.
Each possible distribution of a regenerative composition structure
is thereby described in terms of the drift coefficient $\drift$ and
L{\'e}vy measure $\Om$ of an associated subordinator.
Alternatively, we can transform
$\cal R$ into $\widetilde{\cal R}:= 1 - \exp( - {\cal R } ) \subset [0,1]$ 
and replace the exponential sample by a sample from the uniform distribution
on $[0,1]$.
In this form the construction is an instance
of the {\em ordered paintbox representation} of composition structures, 
developed in \cite{dj91, gnedin97, jp.bmpart}.
\par
Keeping track of only the sizes of parts, and not their order,
every composition structure induces a
{\em partition structure}, that is a sequence of sampling consistent 
{\it partitions} of integers, as studied by Kingman \cite{ki78b, ki80}.
Passing from compositions to partitions is equivalent to passing from the ordered  paintbox  
$\widetilde{\cal R}^c=[0,1]\setminus \widetilde{\cal R}$ 
to {\em Kingman's paintbox} defined by the decreasing sequence of lengths of interval components of 
$\widetilde{\cal R}^c$.
A partition structure is thereby associated with a typically
infinite collection of composition structures, each corresponding to 
a different way of ordering interval components of given lengths.
We show that if one of these composition structures is regenerative,
it is unique in distribution (Corollary \ref{uniquer}). In Section
\ref{structural} we also discuss necessary and sufficient conditions for the 
existence of such a regenerative rearrangement.

\par  Known  examples of regenerative composition structures include the 
compositions associated with the ordered Ewens sampling formula
\cite{dj91},  and those derived from the zero set of a recurrent Bessel process
in \cite{jp.bmpart}. 
The partition structures corresponding to these examples are instances of the two parameter family 
of partition structures studied in \cite{jp.epe,csp}. 
We show in Section \ref{twoparam} that each member of this family, with positive values of parameters,
corresponds to a unique regenerative composition structure.
Also (Theorem \ref{symm}), the only reversible regenerative composition 
structures are the members of this family associated with the interval partition 
of $[0,1]$ generated by excursions of a standard Bessel bridge of dimension 
$2 - 2 \alpha$ for some $\alpha \in [0,1]$.
See also Section \ref{firstex} and \cite{gnedin02three, gnedin03bs}, for further examples of regenerative 
composition structures.

\par Our definition of regenerative composition stuctures is
reminiscent of Kingman's characterisation of the one-parameter Ewens partition 
structure by invariance with respect to deletion of a random part, selected in a size-biased fashion. This property
is called {\em species noninterference}
or {\em neutrality} in the setting of population genetics.
We refer to  \cite{al85, et95, csp}  for 
background on partition structures, exchangeability and related matters.
As shown by James \cite{james03}, another closely related concept, developed 
in the setting of Bayesian nonparametric statistics, is Doksum's 
\cite{doksum74}  notion of a random discrete probability distribution that is 
{\em neutral to the right}.

\par From an algebraic viewpoint, our representation of
regenerative composition structures  is
equivalent to solving
a nonlinear recurrence (Proposition \ref{nonlin}).
The nonlinearity of the recursion reflects the fact that the family of 
probability laws of regenerative compositions is not closed under mixtures.
So unlike the problems of characterising all partition or composition 
structures, the problem of characterising all regenerative composition 
structures is not just a problem of identifying the extreme points of a 
convex set. Still, we show in Section \ref{genrep} that it can be
reduced to such a problem (equivalent to a version of the Hausdorff moment problem), by a suitable 
non-linear transformation. The L{\'e}vy data $(\drift,\Om)$ of the
associated subordinator are thereby encoded in a finite measure on
$[0,1]$.

\section{Compositions and partitions}
This section recalls briefly some background material on composition structures
and their associated partition structures.  See \cite{dj91, gnedin97,  jp.epe, jp.bmpart, csp} for a fuller account.
For a composition structure $({\cal C}_n)$, and a composition
$\lambda=(\ulambda)$ of $n$, define the {\em composition probability function} $p$ by
\begin{equation}
\label{pla} 
p(\lambda) := {\prob} ( {\cal C}_n = \lambda) .
\end{equation}
For each fixed $n$, this function defines a
probability distribution on the set of compositions $\lambda\models n$, 
and these distributions are subject to the following 
 linear relation describing the sampling consistency.
For $\lambda=(n_1,\ldots,n_k)\models n$ and $\mu\models n+1$ we say that $\mu$ {\it extends} 
$\lambda$ and write $\mu\searrow\lambda$
if $\mu$ is obtained from
$\lambda$ by either increasing a part $n_j$ by one or by inserting a part 1 in the sequence $\lambda$.
The sampling consistency amounts to the recursion
\begin{equation}\label{comstr}
p(\lambda)=\sum_{\mu\searrow \lambda}\kappa(\lambda,\mu)p(\mu)\,,\qquad p(1)=1
\end{equation}
where the coefficient $\kappa(\lambda,\mu)$ equals $(n_j+1)/(n+1)$ 
if $\mu$ is obtained by increasing a part $n_j$, and equals $(j+1)/(n+1)$ if $\mu$ is obtained
by inserting a 1 into a row of consecutive ones $1,1,\ldots,1$ of length $j\geq 0$.  

\par 
Regard ${\cal C}_n$ as a way to partition a row
of $n$ identical balls into an ordered series of non-empty boxes,
and independently of ${\cal C}_n$ let the balls be labelled by a uniform
random permutation of the set 
$[n]:= \{1, \ldots, n \}$.
This defines a random  {\em exchangeable ordered partition}
${\cal C}_n^*$ of the set $[n]$
whose distribution is defined as follows. For each {\em particular} ordered partition
of $[n]$ into $k$ classes of sizes $n_1, \ldots, n_k$, say $c^*$,
\eq
\label{pstar}
\prob( {\cal C}_n^* = c^*) = {n \choose n_1, \ldots, n_k }^{-1}
p( n_1, \ldots, n_k)
\en
since the multinomial coefficient is the number of such ordered partitions of
$[n]$, and these are equally likely.
The sampling consistency property of a composition structure $({\cal C}_n)$
means that $({\cal C}_n^*)$ can be constructed {\em consistently}, in the
sense that ${\cal C}_{n-1}^*$ is the restriction of ${\cal C}_n^*$ obtained by deleting element $n$.
Then ${\cal C}_n$ is the ordered record of sizes of classes of ${\cal C}_n^*$, and the
entire sequence $({\cal C}_n^*)$ defines an exchangeable ordered partition of
the set $\mathbb N$ of all positive integers.

\par Ignoring the order of classes yields
a random {\em exchangeable partition} $\Pi$ of the set $\mathbb N$. The restriction 
$\Pi_n$ of $\Pi$ to $[n]$ is obtained by ignoring the order of classes of 
${\cal C}_n^*$.
So for each {\em particular} partition $\pi$ of $[n]$ into $k$ classes whose sizes in some order 
are $n_1, \ldots, n_k$,
\eq
\label{eppf}
\prob( \Pi_n  = \pi) = {n \choose n_1, \ldots, n_k }^{-1} \sum_\sigma p( n_{\sigma(1)}, \ldots, n_{\sigma(k)})
\en
where the sum is over the $k!$ permutations of $[k]$,
corresponding to the $k!$ different ordered partitions $c^*$ of $[n]$ 
derived from the given partition $\pi$ of $[n]$.
This symmetric function of $(n_1, \ldots, n_k)$ is the 
{\em exchangeable partition probability function (EPPF)} of \cite{jp.epe,csp}.
Note by construction that the partition of $n$ defined by the decreasing rearrangement of sizes 
of classes of $\Pi_n$, or of ${\cal C}_n^*$, is identical to  the 
decreasing rearrangement of the parts of ${\cal C}_n$.  Such a sequence of 
random partitions of $n$, subject to a consistency constraint, is called a {\em partition structure}.

\section{Regenerative composition structures}
\label{regener}
Let $({\cal C}_n)$ be a composition structure with composition probability function $p$.
Let $F_n$ denote the size of the first part of ${\cal C}_n$, and denote
the distribution of $F_n$ by
\begin{equation}
\label{qnm} 
q(n:m) := \prob(F_n = m )
= \sum_{(\ulambda)} 1(\nlambda_1=m ) p(\ulambda), \qquad 1 \leq m\leq n,
\end{equation}
where the sum is over all compositions $(\ulambda)$ of $n$,
and $1(\cdots)$ denotes the indicator function which equals 1 if 
$\cdots$ and $0$ else.
We call $q$ the {\it decrement matrix} of the composition structure
$({\cal C}_n)$. 

\begin{proposition}
\label{prop1}
A composition structure $({\cal C}_n)$ 
is regenerative in the sense of Definition
{\em \ref{regendef}} iff for each $n = 1,2, \ldots$ the distribution
of ${\cal C}_n$ is determined by the {\em product formula}
\begin{equation}
\label{produ}
p(\ulambda)= \prod_{j=1}^{\kell} q(\NLambda_j:\nlambda_j)
\end{equation}
for each composition $(\nlambda_1,\ldots ,\nlambda_{\kell})$ of 
$n$, where 
$
\NLambda_j:=\nlambda_j+\cdots +\nlambda_{\kell}
$
and  $q(n : m )$ is the decrement matrix defined by {\em (\ref{qnm})}.
Thus the law of a regenerative composition structure is uniquely determined by 
its decrement matrix.
\end{proposition}
\noindent
\proof
This is easily shown by induction on the number of parts of a composition.
\endpf

\par Note that if $q(2:1)=1$ then $q(n:m)=1(m=1)$,
meaning that each ${\cal C}_n$ is a pure singleton composition, with $p(1,1,\ldots,1)\equiv 1$.
Whereas if $q(2:2)=1$ then $q(n:m)=1(m=n)$ meaning that each ${\cal C}_n$ is a trivial one-part composition
with $p(n)\equiv 1.$ 
These facts are easy to check using (\ref{comstr}) and $p\geq 0$, 
and they are intuitively obvious: $q(2:1)=1$ (respectively $q(2:1)=0$) 
means that two randomly sampled balls never come from the same box (respectively from different boxes).
It can be shown that $q(4:2)>0$ implies
$0<q(n:m)<1$  for all $1\leq m\leq n$ and $n>1$ 
and therefore $0<p(\lambda)<1$ for $\lambda\models n > 1$.
In the case $q(4:2)=0$  and $0<q(2:2)<1$  we have $q(n:1)+q(n:n)=1$ for all $n$, hence $p(\lambda)>0$ only for
compositions of the form $\lambda=(n)$ or $\lambda=(1,1,\ldots,1,k)$ with $k\geq 1$.

\par The product formula (\ref{produ}) identifies 
${\cal C}_n$ with the sequence of decrements of a transient 
Markov chain 
$Q_n:= Q_n(0), Q_n(1), \ldots$
with values in $\{0,\ldots,n\}$. This chain has decreasing paths starting 
from  the  state $Q_n(0) = n$, with the terminal state $0$ 
and time-homogeneous triangular transition matrix $(q(n:n-m), 1\leq m\leq n<\infty)$.
In this interpretation the parts of a composition $\nlambda_1,\ldots,\nlambda_{k}$ are the 
magnitudes of jumps of the chain, while
$(\NLambda_1,\ldots,\NLambda_{k})$ is the path of $Q_n$ prior
to absorbtion.
For example, if ${\cal C}_8 =(3,2,1,2)$, the path of $Q_8$ is
$$
(Q_8(0), Q_8(1), \ldots ) = (8,5,3,2,0,0, \ldots).
$$
\newcommand{\Cn} {{\cal C}_n} 
\newcommand{\Cl} {{\cal C}_n^{<} }
\newcommand{\lambdal} {\lambda^{<} }
\newcommand{\Cr} {{\cal C}_n^{>}  }
\newcommand{\lambdar} {\lambda^{>} }

Consider now the joint law of two compositions derived from a regenerative
composition $\Cn$ by a random splitting,  say $\Cn = (\Cl, \Cr)$, where $\Cl$ 
is a composition of $m(\Cl) \in \{1, \ldots, n \}$, and $\Cr$ is the remaining composition of $n - m(\Cl)$, regarded as a trivial sequence with no elements if $m(\Cl) = n$.
Suppose that the number of parts of $\Cl$ is a {\em randomised stopping time}
of the chain $Q_n$, meaning \cite{ps73} that for each $1 \le k \le n$,
given $\Cn$ with at least $k$ parts, the conditional probability
that $\Cl$ has exactly $k$ parts depends only on the first $k$ parts of $\Cn$.
Equivalently, for each
$\lambda=(n_1,\ldots,n_{\ell})\models n$ and 
each $\lambdal =(n_1,\ldots,n_{k})$ for some $1 \le k \le \ell$,
\eq
\label{new1}
\prob ( \Cl = \lambdal \, | \, \Cn = \lambda ) = f_n( \lambdal)
\en
for some function $f_n$ of compositions of $m$ for $1 \le m \le n$.
The strong Markov property of $Q_n$ then implies that
\begin{itemize}
\item[\rm (i)]
the compositions  $\Cl$ and $\Cr$ are conditionally
independent given $m(\Cl)$, and
\item[\rm (ii)] 
for each $1 \le m < n$, given $m(\Cl)  = m$ the remaining
composition $\Cr$ of $n-m$ is distributed like ${\cal C}_{n-m}$.
\end{itemize}

Conversely, we record the following proposition which applies in particular 
to the splitting scheme defined by \re{new1} with
$f_n( n_1,\ldots,n_{k}) = n_k/n$. In terms of balls in boxes, such a 
split is made just to the right of the box containing a ball picked 
uniformly at random.

\begin{proposition}
\label{prp1}
Suppose a composition structure $({\cal C}_n)$ 
admits a random splitting $\Cn = (\Cl, \Cr)$ for each $n$,
such that {\em \re{new1}} holds with $f_n(m) >0$ for all $1 \le m <n$, 
and {\em (ii) } holds. Then $({\cal C}_n)$ is regenerative.
\end{proposition}
\proof
Let $p$ denote the composition probability function of $(\Cn)$, 
as in \re{pla}. By definition, $(\Cn)$ is regenerative iff for all $1 \le m < n$ and
all compositions $\lambdar$ of $n-m$ 
\eq
\label{pml}
p(m,\lambdar ) = q(n:m) p(\lambdar)
\en
for some matrix $q(n:m)$, which is then the decrement matrix of $(\Cn)$.
Whereas (ii) holds iff for all $1 \le m < n$ and all
compositions $\lambdal$ of $m$ and $\lambdar$ of $n-m$
\eq
\label{pml1}
\sum_ {\lambdal \models m } 
f_n( \lambdal) p( \lambdal, \lambdar) 
= \hat{q}(n:m) p(\lambdar)
\en
for some matrix $\hat{q}(n:m)$, in which case $\hat{q}(n:m) = \prob( m(\Cl ) = m )$.
Assuming that \re{pml1} holds, \re{pml} is obvious for
$m = 1$ with $q(n:1) = \hat{q}(n:1)/f_n(1)$. Proceeding by induction
on $m$, suppose that \re{pml1} holds for all $1 \le m < n$, and that 
\re{pml} has been established with $m'$ instead of $m$ for all
$1 \le m' < m  < n$. Apart from the term $f_n(m) p(m, \lambdar)$,
all terms of the sum in \re{pml1} involve compositions $\lambdal$
all of whose parts are smaller than $m$. So the inductive hypothesis
allows us to write these terms as
$f_n( \lambdal) h_n( \lambdal) p(\lambdar) $ where $h_n(\lambdal )$ is a product 
of entries of the decrement matrix $q$.
Now rearrange \re{pml1} to isolate
the term $f_n(m) p(m, \lambdar)$ on the left, and observe that 
$p(\lambdar)$ is a common factor on the right,
to complete the induction.
\endpf

\par Our aim now is to describe as explicitly as possible all matrices $q$ which 
define a composition structure by means of (\ref{produ}). We start with an algebraic description:
\vskip0.5cm
\begin{proposition} 
\label{nonlin}
A non-negative matrix $q$ is the decrement matrix of some regenerative 
 composition structure iff
$q(1:1)=1$ and 
\begin{eqnarray} \label{q}
q(n:m) =
\frac{m+1}{n+1}\,q(n+1:m+1)+
\frac{n+1-m}{n+1}\, q(n+1:m)+
 \frac{1}{n+1}\,q(n+1:1)\,q(n:m)
\end{eqnarray}
for $1\leq m\leq n$. 
\end{proposition}
\noindent
{\it Proof.} We will show first that the condition (\ref{q}) is sufficient, that is (\ref{q}) and (\ref{produ})
imply (\ref{comstr}). Indeed, assuming (\ref{q}) and (\ref{produ}) 
$$q(n:n)=q(n+1:n+1)+{1\over n+1}q(n+1:n)+{1\over n+1}q(n+1:1)q(n:n)$$ 
implies readily
$$p(n)=p(n+1)+{1\over n+1}p(n,1)+ {1\over n+1}p(1,n)$$
which means (\ref{comstr}) for all one-part compositions. Now suppose (\ref{comstr}) holds for all 
compositions with less than  $k$ parts, and let 
$\lambda\models n$ be a composition with  $k$ parts. Write $\lambda$ in the form $\lambda=(m,\lambda')$ where
$\lambda'\models n-m$. We have by the induction hypothesis and (\ref{produ})
\begin{eqnarray*}
\sum_{\mu\searrow \lambda}\kappa(\lambda,\mu)p(\mu)=
{1\over n+1}p(1,\lambda)+{m+1\over n+1}p(m+1,\lambda')+ 
{n-m+1\over n+1}\sum_{\mu'\searrow \lambda'}\kappa(\lambda',\mu')p(m,\mu')=\\
{1\over n+1}q(n+1:1)q(n:m) p(\lambda')+{m+1\over n+1}q(n+1:m+1)p(\lambda')+ 
{n-m+1\over n+1}q(n+1:m) p(\lambda')
\end{eqnarray*}
which by (\ref{q}) and (\ref{produ}) is equal to $q(n:m)p(\lambda')=p(\lambda)$ and the
induction step is completed.

\par Conversely, assuming (\ref{comstr}) and (\ref{produ}) the recursion (\ref{q}) follows by a similar
argument with   $k=2$.
\endpf

\section{First examples}
\label{firstex}

\paragraph{Example 1} (Geometric sampling \cite{bruss90geom,kar67urn}).
Imagine infinitely many players labeled $1,2,\ldots$ who flip repeatedly the same coin
 with fixed probability $x\in \,]0,1]$ for tails.
In the first round, each of the  players tosses the coin and those who flip tails drop out.
In the second round each of the remaining players must toss again and those who flip tails drop out, and so on.
If we restrict consideration to players labeled $1,\ldots,n$, a composition ${\cal C}_n$ 
arises by arranging the players into 
groups as they drop out. 
These compositions are sampling consistent by exchangeability among the players
and they form a regenerative composition structure because `all rounds are the same'.
Equivalently, we could attribute to each player $j$ an individual value $\xi_j$,
the number of rounds the player remains in the game, and tie the players into blocks by equality of their individual values. The $\xi_j\,$ are independent with same geometric distribution.
The probability that of $n$ players exactly $m$ tie for the minimum value $\min(\xi_1,\ldots,\xi_n)$ is equal to
$$q(n:m)=\frac{{n\choose m} x^m (1-x)^{n-m}}{1-(1-x)^{n}}\,, \qquad m=1,\ldots,n$$ 
which is the binomial distribution conditioned on a positive value. Note that the one-part or 
the pure singleton compositions appear for $x=1$ or $x\downarrow 0$, respectively.
\vskip0.5cm

\par It is the memoryless property which makes the
 geometric distribution work, and sampling from any other {\it fixed} distribution on integers 
would not produce a regenerative 
composition. Still, it is possible to preserve the regenerative feature by randomising the distribution
in a very special way.

\vskip0.5cm
\paragraph{Example 2} (Stick-breaking compositions \cite{dj91, gnedin02three, gnedin03bs,  ho87, james03,  sh85}).
Let $(X_k)$ be independent copies of some random variable
$X$ with $0 < X \le 1$.
Think of $X_k$ as the probability of tails for the $k$th coin.
Modify the algorithm in the previous example by requiring that  
at round $k$ each of the remaining players
must toss the $k$th coin. It is easily seen that
the resulting composition structure is regenerative.
Fixing a group of $n$ players and
conditioning on the number of players that drop out  at the first coin-tossing trial
we obtain the recurrence 
$$q(n:m)={n\choose m}{\ex \,}\left( X^m (1-X)^{n-m}\right) +{\ex \,}(1-X)^n\,q(n:m)$$
resulting in the decrement matrix
\begin{equation}\label{stbr}
q(n:m)=\frac{{n\choose m}{\ex \,}\left(X^m (1-X)^{n-m}\right)}{{\ex \,}\left(1-(1-X)^n\right)} 
\qquad m=1,\ldots,n
\end{equation}
which says that $q(n:\cdot)$ is a mixture of binomial distributions conditioned on a positive value.

\par For example, if $X$ is uniform on $[0,1]$, then
$q(n:m)=n^{-1}$, that is a discrete uniform distribution for each $n$.
More generally, if $X$ has a beta distribution with parameters $(1,\theta)$, $\theta>0$,
 the decrement matrix becomes
\begin{equation}\label{ESF-q}
q(n:m)={n\choose m}\frac{[\theta]_{n-m}\, m!}{[\theta+1]_{n-1}\,n}\,\,,
\end{equation}
where 
\eq
\label{rfac}
[\theta]_n :=\theta (\theta+1)\cdots (\theta+n-1)
\en
is a rising factorial.
The corresponding partition structure is well known to be that defined by the 
Ewens sampling formula \cite{et95}.
The individual values of the players are now only conditionally i.i.d., with conditional distribution
$${\prob}(\xi_j=i\,|\,X_1, X_2, \ldots))=(1-X_1)\cdots(1-X_{i-1})X_i\,.$$

\par Additional randomisation allows the same composition structure to
be defined in another way.
Mark the players
by independent uniform [0,1] random variables $(u_j)$, also independent of $(X_k)$. Consider a random partition of $[0,1]$ into 
intervals by
points 
\eq
\label{sbreak}
Y_k= 1-\prod_{i=1}^k (1-X_i), ~k=1,2,\ldots.
\en
The number of intervals is finite if ${\prob}(X=1) >0$
or infinite otherwise. Group together those players whose
individual marks fall
in the same {\it component} $\,]Y_{k-1}, Y_k[\,$, and maintain the order of groups from the left to the right. 
This sequential algorithm of random interval division is often referred to
as {\em stick-breaking} or as a {\em residual allocation model}.
Note that in the stick-breaking case the partition of $[0,1]$ has a first (leftmost) interval, a second interval, and so on.

\vskip0.5cm 
\paragraph{Example 3} (Brownian bridge \cite{jp.bmpart}). Consider the partition of $[0,1]$ by the set of zeros of a Brownian bridge. This set is perfect, i.e. a compact set with no isolated points.
Given a uniform sample $(u_j)$ group together all sample points which fall into same excursion interval.
This defines a composition structure which is regenerative, by a self-similarity property of the set of zeros.
The decrement matrix is described later by formula \re{q-EP} for $\alpha = \theta = 1/2$.
Unlike the stick-breaking case there is no leftmost interval.

\vskip0.5cm

\paragraph{Example 4}  (Brownian motion, meander case \cite{jp.bmpart}). Same as Example 3 but we take 
the set of zeros of a Brownian motion on $[0,1]$. The collection of intervals is not simply ordered,
but there is a definite last (i.e. rightmost) interval, 
known as the {\em meander} interval, whose right endpoint is $1$.
The decrement matrix is described by formula \re{q-EP} for $\alpha = 1/2, \theta = 0$.

\vskip0.5cm

\paragraph{Example 5} (Myriads of singletons).
Fix $\drift >0$ and a distribution of $X$ on $]0,1]\,$. 
Modify the stick-breaking partition of Example 2 by assuming two types of independent residual allocations. 
At each odd step the stick is broken with residual measure beta$(1,\drift^{-1})$, and at each even step
the stick is broken according to $X$.
That is, consider independent random variables
$Z_1,X_1,Z_2,X_2,\ldots$ with 
$Z_i \stackrel{d}{=}$beta$(1,\drift^{-1})$ and 
$X_i \stackrel{d}{=} X$,
and define
$$
Y_{2k+1}=1-(1- Z_{k+1})\prod_{j=1}^k (1-Z_j)(1-X_j)\,,\qquad
Y_{2k}= 1- (1-X_k)(1-Z_{k})\prod_{j=1}^{k-1} (1-Z_j)(1-X_j).
$$
Consider a random closed set $\widetilde{{\cal R}}$ which includes endpoints $Y_0:=0$ and $1$ and the union of intervals
$[Y_{2k},Y_{2k+1}],$ $k=0,1,\ldots$.
If ${\prob} (X= 1) = 0$ the interval partition has infinitely many components.
\par Draw an independent sample of uniform points $(u_j)$ and define a composition by requiring that 
the sample points which hit components  $[Y_{2k},Y_{2k+1}]$ of $\widetilde{{\cal R}}$ become singletons, while all
those which fall in a particular gap $\,]Y_{2k+1},Y_{2k+2}[\,$ are grouped together.
For $n$ large, a typical composition of $n$  will start with a {\it myriad} of singleton parts
$1,1,\ldots,1$ whose number is of the order of $n$, followed by one part whose
size is of the order of $n$, followed by a myriad, etc.
\par
For $m>1$ conditioning on the number of sample points out of $n$ which fall into $\,]Y_1,Y_2[\,$ leads to a recursion 
$$q(n:m)={n\choose m}{\ex \,}\left((1-Z)^n X^m (1-X)^{n-m}\right) +{\ex \,}\left((1-Z)^n(1-X)^n\right)\,q(n:m)$$
which implies $q$ as in (\ref{stbr}) but with additional term $n\drift$ in the denominator.

\par The total asymptotic frequency of myriads, say $f$, is equal to the Lebesgue measure of $\widetilde{\cal R}$
and satisfies a distributional equation
\eq
\label{frec}
f\stackrel{d}{=} Z_1+(1-Z_1)(1-X_1)f'
\en
where $f',Z_1, X_1$ are independent and $f'\stackrel{d}{=}f$. Analysis of this equation shows that the moments of $f$ are given by a simple formula which 
we record later in  \re{fmoms}.

\section{General representation }
\label{genrep}

\paragraph{Background on subordinators and regenerative sets}

Let $\drift\geq 0$ and $\Om$ be a measure on $\,]0,\infty]$ satisfying
\begin{equation}\label{levcon}
\int_0^{\infty}\min\,(1, z)\,\Om ({\rm d}z)<\infty\,.
\end{equation}
Here and henceforth the integral is over the closed interval 
$[0,\infty]$. There is no mass at $0$ but
we allow the case when $\Om$  gives a positive mass to $z=\infty$.
We also require that either $\drift$ or $\Om$ be nonzero.
Consider a 
Poisson point process on 
$[0,\infty[\,\times [0,\infty]$ with intensity measure Lebesgue$\times \Om$.
Denoting a generic point of the process $(\tau_j,\Delta_j)$,
define the process
\begin{equation}\label{addit}
S_t=\drift\,t +\sum_{\tau_j \le  t}\Delta_j, \qquad t\geq 0 .
\end{equation}
The process $(S_t)$ is a {\em subordinator}, that is a L{\'e}vy process with increasing c{\`a}dl{\`a}g paths,
with $S_0=0$ and $S_t\uparrow \infty$.
For $\rho>0$ let 
$\OW(\rho)$
be the Laplace exponent of the subordinator defined for $\rho \ge 0$ by
$$
{\ex } [ \exp ( - \rho S_t ) ] = \exp [ - t \Phi (\rho )].
$$
Let $\Om( {\rm d} z )$ be the L{\'e}vy measure associated with the subordinator,
and let $\tOm ({\rm d}x )$ be the image of $\Om$ via the transformation 
$x = 1 - e^{-z}$.
According to the L\'evy-Khintchine formula,
\begin{eqnarray}
\label{Ws}
\OW(\rho)
&=&\int_{0}^{\infty} (1-e^{- \rho z})\Om ({\rm d}z)+ \rho \drift \\
\label{denom}
&= & \int_0^1 ( 1 - (1-x)^{\rho})\,\widetilde{\Om}({\rm \,d}x)+ \rho \drift \\
\label{denom1}
& = & \int_0^1 \rho (1-x)^{\rho -1} \tOm [x,1]\,  {\rm d}x   + \rho  \drift .
\end{eqnarray}

Let $${\cal R}=\{S_t, t\geq 0\}^{\rm cl}$$ be the {\em closed range} of the subordinator. 
For a random closed subset $\cal R$ of $[0,\infty]$ let
\eq
G({\cal R},t):= \sup {\cal R} \cap [0,t] \mbox{ and }
D({\cal R},t):= \inf {\cal R} \cap ]t, \infty]
\en
with the usual conventions $\sup \emptyset = 0$ and $\inf \emptyset = \infty$.
Following Maisonneuve \cite{mais83} and 
Bertoin \cite{bert99sub}, 
call $\cal R$ {\em regenerative} if for each $t\in [0, \infty[\,\,$, conditionally on 
$\{D({\cal R},t)< \infty\}$,
 the random set $( {\cal R}- D({\cal R},t) ) \cap [0,\infty]$ is 
distributed like ${\cal R}$ and is independent of 
$[0\,,\,D({\cal R},t)]\cap {\cal R}$. 
The following representation of regenerative sets is fundamental:

\begin{theorem}
\label{maison}
{\em (Maisonneuve \cite{mais83})}
The closed range $\cal R$ of a subordinator $(S_t)$ is a regenerative random subset of
$[0,\infty]$. Moreover, every regenerative random subset ${\cal R}$ of
$[0,\infty]$ has the same distribution as the closed range of some subordinator $(S_t, t \ge 0)$,
whose Laplace exponent $\OW$ is uniquely determined up to constant multiples.
\end{theorem}

\paragraph{Standard exponential sampling}

\par Let $(\epsilon_j)$ be a sequence of independent standard exponential
variables, independent of the subordinator $(S_t)$, 
and let $\epsilon_{1n},\ldots ,\epsilon_{nn}$ be the first $n$ sample points 
$\epsilon_1, \ldots, \epsilon_n$ arranged in increasing order.
Define a partition of the set $\{1, \ldots, n \}$ into blocks of consecutive integers by letting $j$ and $j+1$ belong to different
blocks iff the closed interval $[\epsilon_{jn}\,,\, \epsilon_{j+1,n}]$ 
contains some point of ${\cal R}$, for 
$j<n$. Note in particular that  $\{j\}$ is a singleton block if $\epsilon_{jn}\in {\cal R}$.
Define a composition ${\cal C}_n$ of $n$ by the sequence of counts of block-sizes of this random partition of
$\{1, \ldots, n \}$ into blocks of consecutive integers, from the left to the right.
It is obvious by construction that $({\cal C}_n )$ is a composition structure,
call it the {\em composition structure derived from the subordinator by standard exponential sampling}.

\vskip0.5cm

Introduce the binomial moments 
\begin{eqnarray}
\label{Wnm}
\ow(n:m)&=& {n\choose m} \int_{0}^{\infty} (1-e^{-z})^m\, e^{-(n-m)z}\,\Om ({\rm d}z)+n\drift \,1(m = 1) \\
\label{w-bin-mo}
&=& {n\choose m} \int_0^1 x^m (1-x)^{n-m}\,\widetilde{\Om}({\rm \,d}x)+
 \,n\drift\,1(m =1)  
\end{eqnarray}
for $\tOm ({\rm d}x )$ the image of $\Om( {\rm d} z )$ via $x = 1 - e^{-z}$,
as in \re{Ws}-\re{denom}.
Note by (\ref{levcon}) that the integrals are finite for $1\leq m\leq n$, and that
these quantities are linearly related to the Laplace 
exponent $\OW$ by the elementary identities
\begin{eqnarray}
\label{altbinom}
\OW(n)&=&\sum_{m=1}^n \ow(n:m)\,,~~~~n=1,2,\ldots \\
\label{Phi-nm1}
\ow(n:m) &=&{n\choose m}\sum_{j=0}^m (-1)^{j+1}{m\choose j}\OW(n-m+j)\,, ~~~1\leq m\leq n\,
\end{eqnarray}
where $\Phi(0)=0$.
\vskip0.5cm

\begin{theorem} 
\label{thm1}
\begin{itemize}
\item[\rm (i)]
The composition structure derived from a subordinator by standard exponential
sampling is regenerative, with 
decrement matrix 
\begin{equation}\label{qm}
q(n:m)  = \frac{\ow(n:m)}{\OW(n)} \,.
\end{equation}
\item[\rm (ii)] 
Every regenerative composition structure can be so derived from some subordinator. 
\item[\rm (iii)] 
The L{\'e}vy data $(\drift,\Om)$ of the subordinator is determined uniquely 
up to a positive factor by the regenerative composition structure.
\end{itemize}

\end{theorem}
\vskip0.5cm

To prepare for the proof, we start by recalling some known facts about the 
passage of a subordinator across an independent exponential level. 

\vskip0.5cm
\begin{lemma} {\em \cite{jp.bmpart}}
\label{firstpas}
Let $\epsilon$ be an exponential random variable with rate $\rho$, 
independent of ${\cal R}$ which is the closed range of a 
subordinator $(S_t)$ with Laplace exponent $\OW$.
Let $G_\epsilon:= G({\cal R},\epsilon),$
$D_\epsilon:= D({\cal R},\epsilon),$ and
$\Delta_\epsilon := D_\epsilon - G_\epsilon$, so that
almost surely $\Delta_\epsilon$ is the length of the interval component of
$[0,\infty] \setminus {\cal R}$ which covers $\epsilon$, with
$\Delta_\epsilon = 0$ if $\epsilon \in {\cal R}$.
The random variables $G_\epsilon$ and
$\Delta_\epsilon$ are independent, with  Laplace 
transforms
\eq
\label{gdlap}
{\ex }\exp(-s\, G_\epsilon)=\frac{\OW(\rho)}{\OW(s+\rho)}\,,\qquad 
{\ex }\exp\big(-s\,\Delta_{\epsilon}\big)=\frac{\OW(s+\rho)-\OW(s)}{\OW(\rho)}\,.
\en
\end{lemma}
\vskip0.5cm

\par Note that the second formula in \re{gdlap} is equivalent to
\begin{equation}
\label{fpd}
{\prob}(\Delta_{\epsilon}\in {\rm d}z)=
\frac{(1-e^{- \rho z}) \, \Om({\rm d}z)  + \rho \drift \,\delta_0({\rm d}z)}{\OW(\rho)}\end{equation}
where $\delta_0$ is a unit mass at $0$.
\vskip0.5cm
\noindent
{\it Proof of Theorem {\rm \ref{thm1} (i)}}.
The regenerative property of the composition structure derived from a subordinator 
follows easily from the memoryless property of exponential distribution and the regenerative property of ${\cal R}$ at time 
$D_{1n}:= D({\cal R},\epsilon_{1n})$.
To derive \re{qm},
observe that $\epsilon_{1n}$ is exponential with
rate $n$ and,
by the construction, 
$$q(n:m)=
{\prob}(D_{1n} \,\in\, [\epsilon_{mn}\,,\,\epsilon_{m+1,n}])$$
(with the convention $\epsilon_{n+1,n}=\infty$). 
Let $G_{1n}:= G({\cal R},\epsilon_{1n})$ and $\Delta_{1n}:= D_{1n} - G_{1n}$.
By Lemma \ref{firstpas},
$\Delta_{1n}$ has distribution (\ref{fpd})
for $\rho = n$.
Moreover, given $\Delta_{1n} = z$  with $z >0$, the random variable
$\epsilon_{1n} - G_{1n}$
is distributed like exponential variable $\epsilon(n)$ with rate $n$ conditioned on $\epsilon(n)< z$. 
So the probability that $\epsilon_{1n}$ hits the closed range ${\cal R}$
of the subordinator (causing a singleton) is 
\begin{equation}
\label{singleterm}
{\prob}(D_{1n} =\epsilon_{1n})={\mathbb P}(\Delta_{1n}=0)= \frac{n \drift}{ \OW(n)}
\end{equation}
and given the complementary event that $\epsilon_{1n}$ misses ${\cal R}$,
with $\epsilon_{1n} - G_{1n} =x >0 $ and $\Delta_{1n} =z > x $, 
the conditional probability that $D_{1n} \,\in \,[\epsilon_{mn}\,,\,\epsilon_{m+1,n}]\,$
 equals
$$
{n -1 \choose m-1} ( 1 - e^{-(z-x) })^{m-1} e^{- (z-x) (n-m)}.
$$
So the probability that $\epsilon_{1n}$ finds a gap in ${\cal R}$,
and exactly $m$ of the $n$
exponential variables $\epsilon_1,\ldots,\epsilon_n$ fall in that gap, is
\begin{eqnarray*}
\frac{1}{\OW(n)}
 \int_0^\infty \Om({\rm d}z) \int_0^z n e^{- n x } \, {\rm d}x \, {n -1 \choose m-1} ( 1 - e^{-(z-x) })^{m-1} e^{- (z-x) (n - m)}
\\ \label{}
= \frac{1}{\OW(n)}
 {n \choose m } \int_0^\infty   e^{- (n-m) z } (1 - e^{-z})^m\,\Om({\rm d}z)
\end{eqnarray*}
by application of the formula
$
\int_0^z m e^{-mx } ( 1 - e^{z-x})^{m-1} {\rm d}x = ( 1 - e^{-z})^m
$
which has an immediate interpretation in terms of the order statistics of $m$ independent
exponential variables.
Now (\ref{qm}) follows because $q(n:m)$
 is given by the above formula for $m>1$ and has the additional term $n\drift/\OW(n)$ from (\ref{singleterm}) for $m=1$.
\endpf

To prepare for the proof of the rest of Theorem 
\ref{thm1}, we record a sequence of four preliminary results.
The first is elementary. 

\begin{lemma} 
\label{eltry}
For $1 \le m \le n$ let $\Phi(n:m)$ and $\Phi(n)$ be real variables related 
by {\rm \re{Phi-nm1}} with $\Phi(0) = 0$. Then the identity {\rm\re{altbinom}}
holds.
Moreover, {\rm\re{Phi-nm1}} for $1\leq m\leq n\leq n'$ implies the recursion
\begin{equation}\label{Pascal}
\Phi(n:m)={m+1\over n+1}\,\Phi(n+1:m+1)+{n-m+1\over n+1}\,\Phi(n+1:m), \qquad 1\leq m\leq n<n'\,.
\end{equation}
Conversely, {\rm\re{Pascal}} and {\rm\re{altbinom}} for $1\leq n\leq n'$ 
imply {\rm \re{Phi-nm1}}.
\end{lemma}
A sequence $\Phi$ such that $\ow(n:m)$ defined by \re{Phi-nm1} is non-negative 
for all $n$ and $m$ is known as a {\em completely alternating sequence} 
\cite{berg}, and there is the following integral representation of
such sequences:

\begin{proposition} 
\label{momentsP}
{\em \cite[Proposition 6.12 for $k = 1$, p. 134]{berg} }
A sequence $(\Phi(n), n\geq 0)$ with $\Phi(0 ) = 0$ and $\Phi(n) >0$
for $n >0$ is such that all entries $\Phi(n:m)$ defined by
{\rm (\ref{Phi-nm1})} are non-negative if and only if there is
the integral representation
{\em \re{denom}} for some measure $\tOm$ on $\,]0,1]$ and $\drift\geq 0$. 
Moreover $\tOm$ and $\drift$ are uniquely determined by $\Phi$.
\end{proposition}

\begin{lemma} 
\label{lephi}
Suppose that a sequence of numbers $(\Phi(n), n\geq 0)$ with
$\Phi(0) = 0$ 
satisfies
$\Phi(n) > 0 $ for some $n \leq n'$, and is such that
each entry $\ow(n:m), \,1\leq m\leq n\leq n'\,,$ of the matrix {\rm (\ref{Phi-nm1})}
is non-negative.
Then  $\Phi(n)>0$ for all $1\leq n\leq n'$, and  
the entries of the matrix  
{\em \re{qm}} with $1\leq m\leq n\leq n'$
are non-negative and satisfy {\rm (\ref{q})} for this range of indices.
Moreover,
if the entries $\ow(n:m)$ of the matrix {\rm \,(\ref{Phi-nm1})} are non-negative for arbitrary $n$ then {\em \re{qm}}
is the decrement matrix of some regenerative 
composition structure.
\end{lemma}
{\it Proof.} We apply Lemma \ref{eltry}. 
Dividing (\ref{Pascal}) by $\Phi(n+1)$ and substituting it in the 
to-be-checked (\ref{q}), we transform it by elementary algebra to 
$$\Phi(n+1:1)=(n+1)(\Phi(n+1)-\Phi(n))$$ 
which is true as a special case of (\ref{Phi-nm1}).\endpf

\vskip0.5cm

\begin{lemma} 
\label{leexist}
The decrement matrix of a regenerative composition structure can be represented in the form {\rm(\ref{qm})}, by a matrix 
$(\Phi(n:m),1\leq m\leq n<\infty)$ with non-negative entries satisfying {\rm (\ref{Pascal})} and {\em \re{altbinom}}.
The matrix $\Phi$ is determined by
$q$ uniquely up to a positive factor. 
\end{lemma}
{\it Proof.} The statement is only nontrivial when $0<p(n)<1$ for $n\geq 2$. 
So let us consider a decrement matrix with entries $0<q(n:m)<1$ for $n>1$.
Fix $n'$ and set by definition $\Phi(n':m):=q(n':m)$ for $m=1,\ldots,n'$. Consider the unique
solution $(\Phi(n:m),1\leq m\leq n< n')$ to (\ref{Pascal}) with the  values $q(n':m)$ at level $n'$.
Because $q(n':m)>0$, it is easily seen that
$\Phi(n:m)>0$ for $1\leq m\leq n\leq n'$ and therefore $\Phi(n):=\Phi(n:1)+\cdots +\Phi(n:n)>0$ for $n< n'$
(and $\Phi(n')=1$).
By the first assertion of Lemma \ref{lephi} and the remark before, the elements $\Phi(n:m)/\Phi (n)$ satisfy the recursion (\ref{q}) 
for $n<n'$ and 
for $n=n'$ they coincide with
$q(n':m)$. Thus by the uniqueness of solutions to (\ref{q}) for $n<n'$ with given values at level $n'$ we conclude that
$q(n:m)$ coincides with $\Phi(n:m)/\Phi (n)$ for all $1\leq m\leq n\leq n'$. 

\par Keeping $n'$ fixed, suppose there is another representation $q(n:m)=\widehat{\Phi}(n:m)/\widehat{\Phi}(n),$ $n\leq n'$, 
then $\widehat{\Phi}(n':m)=\widehat{ \Phi}(n') q(n':m)$, thus arguing as above and using linearity we get
 $\widehat{\Phi}(n:m)=\widehat{\Phi}(n')\Phi(n:m)$ for
$1\leq m\leq n\leq n'$. Thus the representation for given $n'$ is unique up to a multiple, and it becomes 
unique subject to a normalisation constraint.

\par Assuming the normalisation $\Phi(1:1)=1$, the finite matrices $(\Phi(m:n), 1\leq m\leq n\leq n')$ constructed 
for each $n'$ are consistent as $n'$ varies, by the uniqueness for each particular $n'$, thus 
they constitute an infinite matrix and the desired representation follows.
\endpf

\vskip0.5cm
\noindent
{\it Proof of Theorem {\rm\ref{thm1} (ii)} and {\rm (iii)}}.
These results follow immediately from Lemma \ref{leexist}  and
Proposition \ref{momentsP}.
\endpf
For an alternate proof of (ii), see 
\cite{prelim_version}.
Also, (iii) can be deduced from Theorem \ref{maison} and a general fact about 
composition structures \cite[Corollary 12]{gnedin97}.

\paragraph{Class frequencies}
If the regenerative composition structure $({\cal C}_n)$ is derived from a subordinator by standard exponential sampling,
the associated composition ${\cal C}^*$ of the infinite set $\mathbb N$ is simply constructed by 
assigning  $i$ and $j$ to different classes
iff the closed interval with endpoints $\epsilon_i$ and $\epsilon_j$  intersects $\cal R$. 
The ordering of classes is maintained according to the order of the 
$\epsilon_j$ associated with the classes.
The random set of positive integers $j$ whose $\epsilon_j$ falls in a particular interval 
component of ${\cal R}^c:=[0,\infty]\setminus {\cal R}$
forms a {\it positive} class, while each $j$ whose $\epsilon_j$ hits ${\cal R}$ forms
a singleton class.
By the law of large numbers, the  probability assigned to an interval component of ${\cal R}^c$
by the standard exponential distribution is the {\it frequency} of the
corresponding class of ${\cal C}^*$, that is the almost sure limit as $n \te \infty$ of
the proportion of elements of $[n]$ which belong to the class.
For instance, if $]a,b[\,\subset {\cal R}^c$  is the interval component which covers $\epsilon_1$,
then for large $n$ the class of ${\cal C}_n^*$ containing element 1 will have
approximately $n(e^{-a} -e^{-b})$ elements, so there will be some part of ${\cal C}_n$ of 
this size.
We note the following Corollary of Theorem \ref{thm1}:

\begin{corollary}
\label{fdist}
Let $f$ denote the random frequency of the union of all singleton
classes in the exchangeable random partition of $\mathbb N$ associated
with a regenerative composition structure with decrement matrix
{\em \re{qm}}. Then
\eq
\label{feq}
f = \drift \, \int_0^\infty \exp({-S_t})\, {\rm d} t
\en
where $(S_t)$ is the associated subordinator with Laplace exponent
$\OW$ and $\drift$ is the drift coefficient of $(S_t)$, and the
distribution of $f$ on $[0,1]$ is determined by the moments
\eq
\label{fmoms}
\ex ( f^n) = { n! \, \drift ^n  \over \mbox{ $\prod_{i = 1}^n$}  \OW (i) }
~~~~~~~~n = 1,2, \ldots .
\en
\end{corollary}
\proof
The derivation from $(S_t)$ by standard exponential sampling gives
$$
f = \int_0^\infty e^{-z} \,\, 1( S_t = z \mbox{ for some } t \ge 0 )\,{\rm d}z
$$
and \re{feq} follows by the change of variable $z = S_t$.
This change of variable follows by noting that 
the function $t\mapsto S_t$ is almost everywhere differentiable with derivative $\drift$.
Formula \re{fmoms} can now
be read from the work of Carmona, Petit and Yor \cite[Proposition 3.3]{cpy97exp},
or derived from \re{frec}.
\endpf

\par Extensive discussion 
of the exponential functional $\int_0^\infty \exp({-S_t})\, {\rm d} t$ is found in \cite{bertyor01,james03}.
See \cite{gpy03regvar} for further applications to regenerative composition structures.

\section{Multiplicatively regenerative sets}
\label{multrep}

\par By mapping $[0,\infty]$ onto $[0,1]$ via $z\mapsto 1-e^{-z}$ we transform
a subordinator $(S_t)$ into  a {\it multiplicative subordinator}
$\widetilde{S}_t:=1-\exp(-S_t)$: for $t'>t$ the ratio $(1-\widetilde{S}_{t'})/(1-\widetilde{S}_t)$ 
has same distribution as $1-\widetilde{S}_{t'-t}$
and is independent of $(S_u, 0 \le u \le t)$.
This construction appears also in \cite{doksum74, gnedin02three, james03}.
The counterpart  of (\ref{addit})  is
$$
\widetilde{S}_t=1-e^{-\drift t} \prod_{\tau_j \le  t} (1-\widetilde{\Delta}_j)
$$
where $\widetilde{\Delta}_j=1-\exp(-\Delta_j)$
and the product is over the atoms 
$(\tau_j,\widetilde{\Delta}_j)$ of a Poisson point process 
in the strip $[0,\infty[\,\times\,[0,1]$, with intensity measure 
Lebesgue$\times\widetilde{\Om}$
where ${\widetilde{\Om}}$ is the image of the measure
$\Om$ via $z\mapsto 1-e^{-z}$.
Note that the mapping preserves order, so that $(\widetilde{S}_t)$ increases from $0$ to $1$. 
\par Let ${\widetilde {\cal R}}:= 1 - \exp( - {\cal R})$ be the closed 
range of the 
multiplicative subordinator
(speaking of closed subsets of $[0,1]$ we shall always mean that
the points $0$ and $1$ are contained in the set).
The transformation $z\mapsto 1-e^{-z}$ takes
an exponential sample $(\epsilon_j)$ into a uniform sample $(u_j)$.
The regenerative composition structure $( {\cal C}_n )$ derived from the subordinator $(S_t)$ by exponential sampling
can now be described as follows: ${\cal C}_n$ is induced by separating the first $n$ uniform 
variables $u_j$ by the points of ${\widetilde {\cal R}}$. 
Note that the frequencies of positive classes derived from
$({\cal C}_n)$ now coincide with the lengths of open interval components of
${\widetilde {\cal R}}^c=[0,1]\setminus
{\widetilde {\cal R}}$, and remaining frequency of singletons $f$,
as in Corollary \ref{fdist}, is the Lebesgue measure of 
${\widetilde {\cal R}}$.

\par For a closed subset $R$ of $[0,1]$ and $z\in \,[0,1[\,$ such that 
$R\,\cap \,\,]z,1[\,\neq \emptyset$, we can define another closed set
\begin{equation}
\label{Rz}
R(z):= \left\{ 
\frac{y-D(R,z)}{1-D(R,z)}
: \, y \in R \cap\, [D(R,z),1]
\right\}
\end{equation}
which is the part of $R$ strictly to the right of $D(R,z)$, scaled back to $[0,1]$. 

\begin{definition}
{\em A random closed set ${\widetilde{\cal R}} \subset [0,1]$ is called {\em multiplicatively regenerative}
if, for each $z\in [0,1[\,\,$, conditionally on $\{D({\widetilde{\cal R}},z)<1\}$
the random set ${\widetilde{\cal R}}(z)$, defined as in (\ref{Rz}),
is independent of $[0\,,\,D({\widetilde{\cal R}},z)]\cap {\widetilde{\cal R}}$, and 
has the same distribution as ${\widetilde{\cal R}}$.}
\end{definition}

The following proposition is easily checked:
\vskip0.5cm
\begin{proposition}
\label{transform}
For random closed sets ${\widetilde {\cal R}} \subset [0,1]$ 
and
${\cal R} \subset [0,\infty]$ related by 
${\widetilde {\cal R}}= 1 - \exp( - {\cal R})$, the random set
${\cal R}$ is regenerative iff
${\widetilde {\cal R}}$ is multiplicatively regenerative.
\end{proposition}

\vskip0.5cm

As a variation of Corollary \ref{crlexp}, a condition for
multiplicative regeneration of a random closed 
subset $\widetilde{{\cal R}}$ of $[0,1]$ can also be given in terms of
a single independent uniform variable.

We associate each composition $(\nlambda_1,\ldots,\nlambda_{\kell})$ of $n$ 
with the finite closed set whose points are  
partial sums of the parts of $\ulambda$ divided by $n$; e.g. the composition
$(4,2,3,1)$ of $10$ is associated with the set $\{0,\,0.4,\,0.6,\,0.9,\,1\}$.
Thus a composition structure $({\cal C}_n)$ is associated with a sequence of random sets $({\widetilde{\cal R}}_n)$.
\noindent
\begin{lemma} 
\label{l59}
{\rm \cite{gnedin97}} Let $({\cal C}_n)$ be a composition structure 
and let $({\widetilde{\cal R}}_n)$ be the associated sequence of random sets.
Then ${\widetilde{\cal R}}_n$ converges almost surely (in the Hausdorff metric)
to some random closed subset ${\widetilde{\cal R}}$,
and $({\cal C}_n)$ is distributed as if by using 
${\widetilde{\cal R}}$ to separate the points in a random sample of uniform $[0,1]$ variables independent
of ${\widetilde{\cal R}}$.
\end{lemma}

From Theorem \ref{thm1}, Proposition \ref{transform} and Lemma \ref{l59} we deduce
\vskip0.5cm
\begin{corollary}
\label{crl5}
The composition structure $({\cal C}_n)$ is regenerative iff $\widetilde{{\cal R}}$ is multiplicatively regenerative. 
\end{corollary}
\vskip0.5cm
As indicated in \cite{prelim_version}, it is also possible to prove
Corollary \ref{crl5} directly, and then retrace the above argument to
obtain an alternate proof of Theorem \ref{thm1}.

\paragraph{A sufficient condition for regeneration}
We note 
that in the usual definition of a regenerative
random subset ${\cal R}$ of $[0,\infty]$, as in Section \ref{genrep},
the independence of the two random 
sets ${\cal R}_t := ( {\cal R}- D({\cal R},t) ) \cap [0,\infty]$ and 
$[0\,,\,D({\cal R},t)]\cap {\cal R}$ for all $t$ can be replaced by
 the apparently weaker condition of independence of the random set
 ${\cal R}_t$ and the random variable $D({\cal R},t)$ for all $t$. This is due to the following result: 

\begin{corollary}
\label{crlexp}
Let ${\cal R}$ be a random closed subset of $[0,\infty]$, let $\epsilon$ be
an exponential random variable with rate $1$ independent of ${\cal R}$, and
let ${\cal R}_\epsilon:= ( {\cal R}- D({\cal R},\epsilon) ) \cap [0,\infty]$.
If ${\cal R}_\epsilon \ed {\cal R}$ and
${\cal R}_\epsilon$ is independent of $D({\cal R},\epsilon)$ then
${\cal R}$ is regenerative.
\end{corollary}
\proof
Let $({\cal C}_n)$ be the composition structure derived from ${\cal R}$
by the standard exponential sampling with variables $(\epsilon_j)$.
Then split $\Cn = (\Cl, \Cr)$, where $\Cl$ is the sequence of non-zero 
numbers of $\epsilon_j$ for $1 \le j \le n$ falling in complementary intervals of 
${\cal R}$ up to and including the count in the interval containing $\epsilon_1$.
This splitting of ${\cal C}_n$ is the example preceding 
Proposition \ref{prp1}, hence
by the assumption on ${\cal R}$ and the memoryless
property of the exponential distribution, it 
satisfies the assumption of Proposition \ref{prp1}.
The conclusion now follows by application of Proposition \ref{prp1}, Theorem \ref{thm1}, 
Corollary \ref{crl5} and Proposition \ref{transform}.
\endpf

\section{Parametrisation of decrement matrices}
\label{parameter}

The representation $q(n:m)=\Phi(n:m)/\Phi(n)$ provides one parametrisation of the regenerative
composition structures in terms of a sequence $(\Phi(n), n\geq 1)$.
To be probabilistically meaningful, this must be the sequence
of evaluations of some Laplace exponent at positive integer values.
But we may also regard the expressions for $q(n:m)$ as a
collection of rational functions in variables $\Phi(n),n \ge 1.$
This section presents some alternative parametrisations of regenerative
composition structures, and discusses their probabilistic and algebraic 
relations to each other.

\subsection{Structural moments}
\label{structural}
One meaningful collection of parameters is the sequence of diagonal entries 
$$p(n)=q(n:n)$$
which starts with $p(1)=1$. 
 We call these diagonal entries of the decrement matrix the
{\it structural moments} of composition structure,
as they coincide with moments of the {\it structural  distribution} $\Sigma$:
$$p(n)= \int_0^1 x^{n-1}\, \Sigma ({\rm d}x)$$
where $\Sigma$ is the distribution of the length 
of the interval component of $\widetilde{\cal R}^c$ containing a given uniform sample point, say $u_1$. 
This random length is the frequency of the class of ${\cal C}^*$ containing element 1, that is a size-biased pick from the  
collection of frequencies \cite{csp}. 
Note from \re{singleterm} and \re{denom}, or from Corollary \ref{fdist},
that the expectation of the total frequency of singletons 
$f = {\rm Lebesgue} ( \widetilde{\cal R} )$ is the measure assigned by
$\Sigma$ to $0$:
$$
\mbox{$\ex (f ) = \Sigma( \{0\}) = \drift/\OW(1) = \drift/\left(\drift+ \int_0^1 x \,\widetilde{\Om}({\rm d}x)\right).$}
$$
\vskip0.5cm
\par  From $p(n)=\Phi(n:n)/\Phi(n),$ 
 by expanding the numerator  by (\ref{Phi-nm1}) we obtain a relation
\begin{equation}\label{Phi-p}
\Phi(n)(p(n)+(-1)^n)=\sum_{j=1}^{n-1} (-1)^{j+1} {n\choose j}\Phi(j)\,,
\end{equation}
which may be seen as a recursion for $\Phi(n),n=1,2,\ldots$. 
Assuming the initial value $\Phi(1)=1$ the recursion has a unique solution,
which is necessarily positive by Lemma \ref{leexist}. Thus the 
recursion (\ref{Phi-p}) allows $q$ to be recovered from $p(n), n=1,2,\ldots$, by first recursively
computing $\Phi(n), n=1,2,\ldots$, then $\Phi(n:m)$ from (\ref{Phi-nm1}) and finally using (\ref{qm}).
Thus we have proved
\vskip0.5cm

\begin{proposition} 
\label{proprat}
{\it A regenerative composition structure is uniquely determined by 
the structural moments $p(n)=q(n:n)$ for $n = 1,2, \ldots$. Each $q(n:m)$ for $1 \le m \le n$
is expressible as a rational function in the variables $p(1)=1,p(2),\ldots,p(n)$.}
\end{proposition} 

\vskip0.5cm
\par To illustrate the result, the first few entries are
\begin{eqnarray*}
q(2:1)&=&1-p(2)\\
q(3:1)&=&\frac{1-3 p(2)+2 p(3)}{1-p(2)}\\
q(3:2)&=&\frac{2 p(2)-3 p(3)+p(2) p(3)}{1-p(2)}\\
q(4:1)&=&\frac{1-5 p(2)+8 p(3)-4 p(2) p(3)-3 p(4)+3 p(2) p(4)}{1-2 p(2)+2 p(3)-p(2) p(3)}\\
q(4:2)&=&\frac{3 p(2)-9 p(3)+6 p(2) p(3)+6 p(4)-9 p(2) p(4)+3 p(3) p(4)}{1-2 p(2)+2 p(3)-p(2) p(3)}\\
q(4:3)&=&\frac{3 p(3)-3 p(2) p(3)-4 p(4)+8 p(2) p(4)-5 p(3) p(4)+p(2) p(3) p(4)}{1-2 p(2)+2 p(3)-p(2) p(3)}\,\,
\end{eqnarray*}
The complexity of such formulas increases rapidly with $n$. 

\par In general,  structural moments  do not determine a composition structure uniquely, 
because they do not even determine the associated partition structure.  
See \cite{csp} for further discussion.
Since uniqueness  does hold in the special case of regenerative composition structures, it is natural 
to seek a characterisation of structural moments in 
this case. 
There is the following immediate consequence of
 Proposition \ref{proprat} and Lemma \ref{lephi}:

\vskip0.5cm

\begin{corollary}
\label{characp}
A sequence $p(n), n=1,2,\ldots$   with $p(1)=1$ and $0<p(n)<1$ for  $n>1$ 
is a sequence of structural moments of some regenerative composition structure
if and only if the following conditions are fulfilled:
\begin{itemize}
\item[{\rm (i)}]  the sequence $\Phi(n), n=1,2,\ldots$ defined by the recursion {\rm (\ref{Phi-p})} with
$\Phi(1)=1$ is positive, and
\item[{\rm (ii)}]  each $\Phi(n:m), 1\leq m\leq n<\infty$ defined by {\rm (\ref{Phi-nm1})} is non-negative. 
\end{itemize}
If this is the case, 
\begin{eqnarray*}
p(n)=\frac{\int_0^1 x^n\, 
\widetilde{\Om}({\rm d}x)}{\int_0^1 (1-(1- x)^n ) \, \widetilde{\Om}({\rm d}x)+n\,\drift}\qquad n>1
\end{eqnarray*}
for some $\drift\geq 0$ and some measure $\widetilde{\Om}$ on $]0,1]$ with finite first moment.
\end{corollary}

\vskip0.5cm
\par {\bf Remark.} We know that $p(n),n=1,2,\ldots$ is a moment sequence  from the general facts about
 partition structures, or from the interpretation of $p(n)$ as the probability that $n$ balls fall in the same box.
From an analytical perspective, 
it does not seem  obvious that the nonlinear tranform given by $p(n)=\Phi(n:n)/\Phi(n), n=1,2,\ldots$ indeed
yields a completely monotonic sequence
for arbitrary Laplace exponent.
\vskip0.5cm  
\par Because the structural moments are determined by the (unordered) partition structure, Proposition \ref{proprat}
and Kingman's representation of partition structures \cite{ki78b} imply:
\vskip0.5cm  
\begin{corollary} 
\label{uniquer}
Each distribution of an infinite exchangeable partition of $\mathbb N$ (which can be identified
with a partition structure)
corresponds to at most one regenerative composition structure.
Equivalently, for each distribution of a decreasing sequence $(Y_j)$ with $Y_j \ge 0$ and $\sum Y_j\leq 1$,
there exists at most one distribution for a multiplicatively regenerative set $\widetilde{\cal R}\subset [0,1]$ such that the ranked lengths of 
interval components of $\widetilde{\cal R}^c$ are distributed like $(Y_j)$.
\end{corollary}

\vskip0.5cm

\par A constructive method to verify if a given exchangeable partition of $\mathbb N$ is induced 
by a regenerative composition structure amounts to computing $q$ from the structural moments,
and then checking that the given EPPF coincides with the EPPF computed by formulas \re{produ} and \re{eppf}.
\par The general problem of characterising  structural distributions  
of partition structures was posed by Pitman and Yor \cite{py95rdd}.
The characterisation of structural distributions of regenerative
composition structures provided by Corollary \ref{characp} leaves open
the following question: given the collection of structural moments of a 
regenerative composition, or given its Laplace exponent $\OW$,
describe in some way how the classes of the associated unordered partition should be arranged to produce the composition.
We answer some restricted forms of this question in the next section, but do not see how to answer it in any generality.

\subsection{Singleton probabilities}

Instead of the event `$n$ balls fall in same box',
consider the event `$n$ balls fall in $n$ different boxes'.
Let $e(n)$ be the probability of this event, that is
$$e(n):=p(1,1,\ldots,1)=q(n:1)q(n-1:1)\cdots q(2:1).$$
By the definition and from the representation (\ref{qm}) we derive
$$
{e(n)\over e(n-1)}=q(n:1)=n\left(
1-{\Phi(n-1)\over \Phi(n)}
\right)
$$
which can be read as
\begin{equation}
\label{eP}
{\Phi(1)\over \Phi(n)}=\prod_{j=2}^n \left( 1-{e(j)\over j\, e(j-1)}\right).
\end{equation}
This shows that any one of the sequences $(e(n),n>0)$, $(q(n:1),n>0)$ or  $(\Phi(n)/\Phi(1),n>0)$ 
uniquely determines each of the other two sequences.

\par As is seen from (\ref{Phi-nm1}) and (\ref{eP}), in the variables $q(n:1),n=1,2,\ldots$ the elements of decrement matrix become polynomials
\begin{equation}\label{q-q1}
q(n:m)={n\choose m}\sum_{j=0}^m (-1)^{m-j+1}{m\choose j}\prod_{k=0}^{j-1}\left(1-{q(n-k:1)\over n-k}\right)\,,
\end{equation}
to be compared with the rational functions of structural moments considered
in Subsection \ref{structural}.
For example
$$q(4:2)=2\,q(3:1)-\frac{3}{2}\,q(4,1)-\frac{1}{2}\,q(3:1)\,q(4:1).$$

\par The definition of $e(n)$ makes sense for a general partition structure. Thus to check if a given partition
structure is induced by a regenerative composition structure, we can use the above formulas
to translate $e(n),n>0,$ into $q$ and then compare the EPPF resulting from (\ref{produ}), (\ref{eppf}) with the
given EPPF. 
In particular, if 
a regenerative rearrangement is possible, the sequences $(p(n), n>0)$ and $(e(n), n>0)$
must be computable from each other, as appears by eliminating the variables $\Phi$ from $p(n)=\Phi(n:n)/\Phi(n)$ and
  (\ref{eP}).

\vskip0.5cm

\section{The two-parameter family}
\label{twoparam}
\subsection{General setup}
Consider the
{\em $(\alpha,\theta)$-partition structure} determined by following formula 
of \cite{jp.epe,csp}
for the
distribution of $\Pi_n$, an exchangeable partition of $[n]$:
for each particular partition $\pi$
of $[n]$  into $k$ classes of sizes $n_1, \ldots, n_k$
\eq
\label{peppf}
\prob( \Pi_n = \pi ) = {\prod_{ i = 1}^{k-1} ( \theta + \alpha i)  \over [1 + \theta]_{n-1} } \prod_{i = 1}^k [ 1 - \alpha ]_{n_i - 1} 
\en
where the notation \re{rfac} is used for rising factorials.
This formula defines a partition structure for $0\leq \alpha<1$ and $\theta\geq 0$, and also for 
some $(\alpha,\theta)$ with either $\alpha <0$ or $\theta < 0$. We wish to 
establish if this partition structure can be associated with some regenerative composition structure.

\par Following the method in Section \ref{parameter} we first compute $e(n)$ as a special case of (\ref{peppf}):
$$e(n) =p(1,1,\ldots,1)=\prod_{j=0}^{n-1}{\theta+\alpha j\over \theta+j}$$
which leads by application of (\ref{eP}) to 
\begin{equation}\label{P-ta}
{\Phi(n)\over \Phi(1)}={n[\theta+1]_{n-1}\over [2+\theta-\alpha]_{n-1}}\,.
\end{equation}
This yields,
by virtue of (\ref{Pascal}) or (\ref{Phi-nm1}),  the formula
$${\Phi(n:m)\over\Phi(1)}={n\choose m}{[1-\alpha]_{m-1}\over   [2+\theta-\alpha]_{n-1}}\,{   [\theta+1]_{n-1}\over
 [\theta+n-m]_m}\,\left((n-m)\alpha +m\theta\right).$$
Therefore
\begin{equation}\label{q-EP}
q(n:m) ={ \ow(n:m) \over \OW(n)}
= 
{ n \choose m } {[1-\alpha]_{m-1} \over [ \theta + n - m]_m} { ( (n - m ) \alpha + m \theta ) \over n }\,\,.
\end{equation}
Since $q$ in (\ref{q-EP}) is non-negative exactly when
$0\leq \alpha<1$ and $\theta\geq 0$ we conclude that 
$q$ is the decrement matrix of a
regenerative composition structure for precisely this range of parameters.

\par Observe that the resulting formula 
\begin{equation}\label{p-EP}
p(n)=q(n:n) = {[1-\alpha]_{n-1}\over [1 + \theta]_{n-1}}
\end{equation}
yields the moments of beta$(1 - \alpha, \alpha + \theta)$, which is  the structural 
distribution for {\it all} members of the two-parameter family of partition structures.

\par 
Adopting the normalisation $\Phi(1)={\rm B}(1-\alpha, 1+\theta)$,
where
$$
{\rm B}(a,b) := \Gamma(a) \Gamma(b) / \Gamma(a+ b)
$$
the Laplace exponent extending  
(\ref{P-ta}) becomes
\eq
\label{althlap}
\OW (s) = s {\rm B}(1-\alpha, s + \theta)  .
\en
The corresponding measure is determined by the formula
\eq
\label{althmeas}
\tOm [x,1] = x^{-\alpha} (1-x)^{\theta} , \qquad 0 < x < 1.
\en
It remains to check that the partition structure induced by
this regenerative composition structure is  given by (\ref{peppf}). This is done in the following theorem:
\vskip0.5cm

\begin{theorem}
\label{alth}
For $0 \le \alpha < 1$ and $\theta \ge 0$ the distribution of the exchangeable random partition $\Pi_n$ of $[n]$ derived from 
the regenerative composition structure with Laplace exponent {\em \re{althlap}} is that of an $(\alpha,\theta)$
partition defined by formula {\em (\ref{peppf})}. For other values of $(\alpha, \theta)$, besides the
limiting case $(1,\theta)$ for $\theta \ge 0$ which generates the pure 
singleton partition,
there is no regenerative composition structure which generates an $(\alpha, \theta)$-partition structure.
\end{theorem}

\proof
By the above discussion we can restrict  consideration to the case $0 \le \alpha < 1$ and $\theta \ge 0$.
By application of formulas \re{eppf}, \re{produ} and \re{q-EP}, the EPPF derived from the
regenerative composition structure with Laplace exponent \re{althlap} is a sum of $k!$ terms of the
form
$$
{1 \over [\theta]_{n} } \prod_{i = 1 }^k [ 1 - \alpha]_{n_i - 1} { (N_i - n_i) \alpha + n_i \theta  \over N_i }
$$
where the sequence $(n_1, \ldots, n_k)$ and its tail sums 
$N_i = \sum_{j = i}^k n_j$ must be replaced by permutations of the sequence 
and correspondingly transformed tail sums.
To match up with \re{peppf}, it just has to be checked that the corrresponding
sum of $k!$ terms derived from
\eq
\label{fack}
\prod_{i = 1 }^k { (N_i - n_i) \alpha + n_i \theta  \over N_i ( (k-i) \alpha + \theta ) }
\en
equals 1. But this is easily verified together with the probabilistic
interpretation given in the following corollary.
\endpf

\begin{corollary}
\label{crlord}
In the setting of the previous theorem, given that the blocks of $\Pi_n$ are of sizes 
$n_1, \ldots, n_k$ when put in some arbitrary order, and given that
the first $i-1$ of these blocks are the first $i-1$ blocks of 
the ordered partition ${\cal C}_n^*$, the conditional probability that this 
coincidence continues for one more step is the $i$th factor in {\em \re{fack}}.
\end{corollary}
Put another way, given block sizes $n_1, \ldots, n_k$
and that the first $i-1$ blocks have been picked to leave blocks of
sizes $n_j$ for $i \le j \le k$, the next block is the block of index $j$ with
probability proportional to $(N_i - n_j) \alpha + n_j \theta$.

\vskip0.5cm

Several particular instances of the above results are known, as indicated in the following discussion
of special cases.

\subsection{Case $(0,\theta)$ for $\theta \ge 0$} 
In this case the measure $\tOm$ in
\re{althmeas} is a probability measure, the beta$(1,\theta)$ distribution. 
So the above theorem and its corollary reduce 
to the well known fact that the ordered Ewens formula associated with beta$(1,\theta)$ stick-breaking 
puts its parts in a size-biased random order \cite{dj91}.

\subsection {Case $(\alpha,0)$ for $0 < \alpha <1$}
\label{alpha-zero}
In this case 
$$\tOm ({\rm d}x)=   \alpha x^{-\alpha -1} {\rm d}x+ \delta_1 ({\rm d}x)$$
is a measure with a beta density on $\,]0,1[\,$ and a unit atom at $1$.
The product formula \re{produ} reduces to 
$$
p(\ulambda)=\nlambda_{\kell} \alpha^{\kell-1}\prod_{j=1}^{\kell}{[1-\alpha]_{\nlambda_j-1}\over \nlambda_j!},
$$
which is identical to the formula in \cite[Equation (28)]{jp.bmpart}.
By comparision of these two formulas, the random composition is in this case is
identical in distribution to that generated by ${\cal R}_\alpha \cap [0,1]$ 
where ${\cal R}_\alpha $ is the range of a stable subordinator of index 
$\alpha$. 
In particular, ${\cal R}_\alpha$ can be realised as the zero set of a 
Bessel process of dimension $2 - 2 \alpha$.
For $\alpha = 1/2$ this is the zero set of a standard Brownian motion.
\par 
The decrement matrix $q$ in this case has the special property that
there is a probability distribution $f$ on the positive integers such that 
\eq
\label{qnf}
q(n:m) = f(m) \mbox{ if } m < n \mbox{ and } q(n:n) = 1 - \sum_{m = 1}^{n-1}f(m).
\en
Specifically,
\begin{equation}\label{f}
f(m)= { \alpha [ 1 - \alpha]_{m-1} \over m! }
\end{equation}
and hence $q(n:n)= [ 1 - \alpha]_{n-1} /(n-1)!$.
The work of Young \cite{young95phd} shows
that the only non-degenerate regenerative composition structures  with a decrement matrix of the
form \re{qnf}, for some probability distribution $f$ on the positive
integers, are those with $f$ of the form \re{f}, obtained
by uniform sampling from
${\cal R}_\alpha \cap [0,1]$ for some $0 < \alpha < 1$.

\par The multiplicative regeneration property of
${\cal R}_\alpha \cap [0,1]$ is an immediate consequence of the standard regeneration
and self-similarity properties of ${\cal R}_\alpha$ as a subset of $[0,\infty]$.
It implies that ${\cal R}_\alpha \cap [0,1]$ has the same distribution 
as the closure of $\{ 1 - \exp(-S_t), t \ge 0 \}$ where $(S_t)$ is
a subordinator with no drift and L\'evy measure 
$$\Om({\rm d}z)=\alpha (1-e^{-z})^{-\alpha-1}e^{-z}{\rm d}z+ \delta_\infty({\rm d}z)$$
on $[0,\infty]$ which is
the image of $\tOm$ via $x \mapsto - \log (1-x)$, so $\Om$ has
an atom of mass $1$ at $\infty$.
\par As a check, let $\tau := \inf \{t : S_t = \infty \}$,
which is the exponential
time with rate $1$ when the subordinator jumps to $\infty$. 
Then, by application of the transformation and the L\'evy-Khintchine formula, if we let
$G := \sup {\cal R}_\alpha \cap [0,1[\,$, then we find for $s >0$
$$
{\ex } (1 - G)^s = {\ex \,} ( \exp(- s S_{\tau-})) = {1 \over \OW(s)} =
{ {\rm B}( 1 - \alpha + s, \alpha ) \over {\rm B}( 1 - \alpha, \alpha ) }.
$$
This confirms the well known fact that the distribution of
$1 - G$ is beta$(1-\alpha,\alpha)$.
It may also be observed, using properties of the local time process 
$(L_t, t \ge 0)$ associated with ${\cal R}_\alpha$, as discussed in 
\cite{neretin96}, that the exponential time $\tau$ can be represented as
$$
\tau = c_\alpha \int_0^1 (1-t)^{- \alpha} {\rm d} L_t
$$
for some constant $c_\alpha$ depending on the normalisation of the local time process.
The fact that this local time integral has an exponential distribution was derived by an
analytic argument in \cite[Corollary 3.4]{MR1702241}.

As discussed in \cite{jp.bmpart}, the length of the last interval component $\,]G,1[\,$
of the complement to ${\cal R}_\alpha \cap[0,1]$ is a size-biased pick from the collection of the interval lengths,
and conditionally on $G$ the remaining interval components are in symmetric order; moreover these properties are inherited by the compositions 
of $n$ for every $n$.
Corollary \ref{crlord} in this case is new. It makes precise another sense in which given the partition of $n$ generated
by ${\cal R}_\alpha \cap [0,1[\,$, the smaller blocks tend to come first in the composition of $n$.

\subsection{Case $(\alpha,\alpha)$ for $0 < \alpha <1$}
\label{alal}
Passing to the variable $z=-\log (1-x)$ we see from \re{althmeas} that the associated
regenerative subset  of $[0,\infty]$ has zero drift and L{\'e}vy measure
$$
\Om ({\rm d} z)
= \alpha (1-e^{-z})^{-\alpha-1} e^{-\alpha z}{\rm\, d}z\qquad z\in\, ]0,\infty[\,.
$$
It can be read from \cite{py96ou} that such a regenerative set is generated
as the zero set of 
the squared Ornstein-Uhlenbeck process $(X_t)$ of dimension $2 - 2 \alpha$ driven 
by the stochastic differential equation ${\rm d}X_t = 2 \sqrt{X_t}\, {\rm d} B_t + (2 - 2 \alpha - X_t){\rm d}t$ 
where $(B_t)$ is a standard Brownian motion, and that the image of 
this regenerative set via $x = 1 - e^{-z}$ is the zero set of a
Bessel bridge of dimension $2 - 2 \alpha$. 
In case $\alpha=1/2$ this is a Brownian bridge, as in Example 3.
In the notation introduced in the discussion of the previous case,
this corresponds to conditioning ${\cal R}_\alpha \cap[0,1]$ on the event
$1 \in {\cal R}_\alpha$. This can be rigorously understood by
first conditioning on $ G \in [1-\epsilon ,1 ]$ and then taking a weak limit 
as $\epsilon \downarrow 0$. 
\par The decrement matrix in this case has the special property that 
\begin{equation}\label{renew}
q(n:m)=\frac{f(m)\,r(n-m)}{r(n)}\,
\end{equation}
where $f$ is given by (\ref{f}) and 
$r(n)=[\alpha]_n/n!$
is the probability that a random walk on positive integers with step distribution
$f$ visits $n$. Equivalently, the composition probability function is
\begin{equation}\label{renew1}
p\,(\ulambda)=\frac{\prod_{j=1}^k f(\nlambda_j)}{r(n)} 
\end{equation}
or more explicitly
\begin{equation}\label{alpha-alpha}
p\,(\ulambda)={n! \over [\alpha]_n} \alpha^k \prod_{j=1}^k  { [1 - \alpha]_{n_i - 1} \over n_i ! } .
\end{equation}
It follows from a result of Kerov \cite{kerov95} that 
the decrement matrix of a non-degenerate regenerative composition
structure can be expressed in the form {\rm(\ref{renew})} for some functions $f$ and $r$ iff it is of the form \re{alpha-alpha} for some $\alpha \in ]0,1[$.
The same conclusion is also a consequence Theorem \ref{symm} in the next section.
The conclusion of Corollary \ref{crlord} in this case is that
given the partition of $[n]$ the block sizes appear in ${\cal C}_n$ in a uniform random order.
This can be seen directly from the symmetry of formula \re{renew1} as 
a function of $(\ulambda)$.

\subsection{Case $(\alpha,\theta)$ for $0<\alpha<1,\, \theta>0$}

It is known \cite{py95rdd,py95pd2,jp97cmc} that an $(\alpha,\theta)$ partition of $\nints$ can be constructed as follows.
First construct a $(0,\theta)$ partition of $\nints$, then shatter each class of this partition
according to an independent $(\alpha,0)$ partition.
This operation restricts naturally to $[n]$ for each $n$, and can be interpreted in terms of
a fragmentation operation on the frequencies of classes.  
This result can be lifted to the level of regenerative composition structures as follows.

\begin{theorem}
\label{thalt}
For $0 < \alpha < 1$ and $\theta >0$, let $Y_0 = 0$ and let $0 < Y_1 < Y_2 <  \ldots$ be defined by the independent
stick-breaking scheme {\em \re{sbreak}} for $X$ with beta$(1,\theta)$ distribution,
let ${\cal R}_{\alpha}(i)$ for $i = 1,2,\ldots$ be a sequence of
independent copies of the range ${\cal R}_{\alpha}$ of a stable subordinator, and define
a random closed subset $\widetilde{{\cal R}}_{(\alpha,\theta)}$ of $[0,1]$ by
$$
\widetilde{{\cal R}}_{(\alpha,\theta)} = \{1\} \cup \bigcup_{i = 1}^\infty \left( [Y_{i-1},Y_i] \cap [Y_{i-1} + {\cal R}_{\alpha}(i)] \right) .
$$
Then $\widetilde{{\cal R}}_{(\alpha,\theta)}$ is a multiplicatively regenerative random subset of $[0,1]$, which can
be represented as $\widetilde{{\cal R}}_{(\alpha,\theta)} = 1 - \exp(-{\cal R}_{(\alpha,\theta)})$ where
${\cal R}_{(\alpha,\theta)}$ is the range of a subordinator with Laplace exponent {\em \re{althlap}},
and the composition structure obtained by uniform random sampling from $\widetilde{\cal R}_{(\alpha,\theta)}$ is regenerative
with decrement matrix {\em \re{q-EP}}.
\end{theorem}
\proof
It is easily checked, using the muliplicative regeneration of the stick-breaking scheme, and the self-similarity of
${\cal R}_{\alpha}$, that $\widetilde{{\cal R}}_{(\alpha,\theta)}$ is multiplicatively regenerative.
The description of the Laplace exponent then follows from Proposition \ref{proprat}, since the structural distribution is easily identified.
\endpf

The particular case $\alpha = \theta$ of Theorem \ref{thalt} is largely contained in the work of Aldous and Pitman \cite{ap01d}.
In particular, for $\alpha = \theta = 1/2$ this construction of the zero set of a Brownian bridge  plays a key role
in the asymptotic theory of random mappings developed in \cite{ap92} and \cite{ap01d}.

\section{The Green matrix}
\label{green}
For a given composition probability function \re{pla},
the {\it Green matrix} is defined by the formula 
  $$g(n,j)=\sum_{\lambda\models n, \,j\in \{N_i\}} p(\lambda),\qquad 1\leq j\leq n<\infty$$
where the summation is over all compositions  $\lambda=(n_1,\ldots,n_k)\models n$ which have integer $j$ among tail sums
$N_j=n-n_1-\ldots -n_{j-1}$ (where we set $n_0=0$).
Recalling the interpretation of a regenerative composition structure 
as a consistent family of Markov chains $Q_n, n=1,2,\ldots$,   as in Section \ref{regener},
$g(n,j)$ is the chance that $Q_n$ with transition matrix 
$q$ and initial state $n$ ever visits
state $j$. In particular, $g(n,n)=1$.
\vskip0.5cm
\noindent
{\bf Example.} For the 2-parameter family we have for $1 \le j \le n$
\begin{itemize}
\item[{(i)}] for $(0,\theta)$
$$g(n,j)={\theta\over j+\theta}$$  as is well known;
\item[{(ii})] for $(\alpha,0)$
$$g(n,j)={[\alpha]_{n-j}\over (n-j)!}$$ 
which by \re{qnf}-\re{f}
is the probability that a particular random walk with negative 
increments started at level $n$ ever visits state $j$;

\item[{(iii)}] for $(\alpha,\alpha)$
$$g(n,j)={{n\choose j}[\alpha]_{j}\over (\alpha+n-j)\cdots (\alpha+n-1)}$$
which is the probability of the same event for the random walk of the
previous case conditioned to hit $0$. 
\end{itemize}

\vskip0.5cm
\begin{lemma} The Green matrix of a regenerative composition structure is the unique solution of the recursion
\begin{equation}\label{rec-g}
g(n,j)={j+1-q(j+1:1)\over n+1}\,g(n+1,j+1)+{n+1-j\over n+1}\,g(n+1,j)
\end{equation}
with boundary condition $g(n,n)=1$. 
\end{lemma}
{\it Proof.} The path of the chain $Q_n$, defining a composition
of $n$,  is 
obtained via random deletion of a state from $1,2,\ldots,n+1$, then restricting a path of $Q_{n+1}$ to the undeleted
 states
and re-labeling the states by ranking them from $1$ to $n$. The event `$Q_n$ visits $j$' occurs
when either $Q_{n+1}$ visits $j$ and one of the states $j+1,\ldots,n+1$ is deleted
(in which case state $j$ retains the label) or $Q_{n+1}$ visits $j+1$ and one of the states $1,\ldots,j+1$ is 
deleted (if
state $j+1$ is not deleted it changes the label to $j$). The first event has probability $g(n+1,j)(n+1-j)/(n+1)$
and the second $g(n+1,j+1)(j+1)/(n+1)$. The events are not disjoint and their intersection is the event
`$Q_{n+1}$ visits both $j+1$ and $j$, and state $j+1$ is deleted' which has probability $g(n+1,j+1)q(j+1:1)/(n+1)$.
The uniqueness claim is obvious from the recursion.
\endpf
\vskip0.5cm

\par Next result gives an explicit formula for the Green matrix in terms of the representation (\ref{qm}) via Laplace
exponent.
\vskip0.5cm
\begin{theorem} The Green matrix of a regenerative composition structure is 
\begin{equation}\label{Green}
g(n,j)=\Phi(j) {n\choose j}\sum_{a=0}^{n-j} {n-j\choose a}{(-1)^a \over \Phi(j+a)}\,.
\end{equation}
\end{theorem}
{\it Proof.} In view of 
$$q(j+1:1)=(j+1)\left(1-{\Phi(j)\over \Phi(j+1)}\right)$$ 
the first factor in the right side of (\ref{rec-g}) equals $(j+1)\Phi(j)/((n+1)\Phi(j+1))$. Substituting 
this and (\ref{Green})
into
(\ref{rec-g}), and canceling the common factor
${n\choose j}\Phi(j)$ the to-be-checked recursion follows from the identity
$$\Delta^{n-j+1} s (j)= \Delta^{n-j}
s (j+1)-
\Delta^{n-j} s (j)$$
where $\Delta$ is the forward difference operator $\Delta s (i) := s(i+1) - s(i)$ and $s$ is the sequence $s(i) = 1/\Phi(i)$ for $i \ge 1$.
\endpf
\vskip0.5cm

\par We give one application of the formula. Let $L_n$ be the {\it last} part of ${\cal C}_n$.
In the event $\{L_n=j\}$ the chain $Q_n$ visits state $j$ and then has the last positive decrement $j$.
The distribution of the last part follows from this observation and 
(\ref{Green}):
\begin{equation}\label{Ln}
{\mathbb P}(L_n=j)=g(n,j)q(j:j)=\Phi(j:j) {n\choose j}\sum_{a=0}^{n-j} {n-j\choose a}{(-1)^a \over \Phi(j+a)}\,.
\end{equation}
In particular, normalising by $\Phi(1)=1$ for simplicity,
\eq
\label{lastone}
\prob ( L_n = 1) = 
n \left[ 1 - \sum_{ k = 2 }^{n } { n - 1 \choose k - 1 } { (-1)^k \over \OW (k) } \right]\,.
\en

\section{Symmetry}

Each composition structure 
$(  {\cal C}_n ) $ has a dual 
$(  \hat{{\cal C}}_n ) $, where
$\hat{{\cal C}}_n $ is the sequence of parts of
$ {\cal C}_n $ in reverse order. If
$ ({\cal C}_n )$ is derived by uniform sampling from a random closed set
$\widetilde{\cal R} \subset [0,1]$, then $\hat{{\cal C}}_n $ is derived similarly from $1 - \widetilde{\cal R}$.
If $(  {\cal C}_n ) $ is regenerative, and 
so is $(\hat{{\cal C}}_n )$, then 
$(  {\cal C}_n ) $ and $(\hat{{\cal C}}_n )$ must 
be identical in distribution, by Corollary \ref{uniquer}.
Equivalently, $\widetilde{\cal R} \ed 1 - \widetilde{\cal R}$, in which
case we call the composition structure {\em reversible}.
Two degenerate examples are provided by $\widetilde{\cal R}  = \{0\} \cup \{1\}$
and $\widetilde{\cal R} = [0,1]$.
The existence of regenerative composition structures which are 
non-degenerate and reversible is quite surprising and counter-intuitive, because the ideas of stick-breaking and multiplicative 
regeneration suggest that typical interval sizes should decay in some
sense from the left to the right. However, it is evident from 
the formula \re{renew1} that for every $0 < \alpha < 1$
the regenerative composition structure
associated with an $(\alpha,\alpha)$ partition
is reversible.
Indeed, this composition structure is {\em symmetric}, meaning that
the composition probability function is a symmetric function of $(\ulambda)$
with respect to all permutations of the arguments, for each $k$. The equivalent
condition on $\widetilde{\cal R}$ is that the interval components of
the complement of $\widetilde{\cal R}$ form an exchangeable interval partition of $[0,1]$,
as defined in \cite{al85}.
We note in passing that a large family of symmetric composition structures 
was derived from the jumps of a subordinator in \cite{pitman02pk}.
See also \cite{gnedin98pd}.

\begin{theorem}
\label{symm}
Let $(  {\cal C}_n )$ be the regenerative composition structure derived by uniform sampling from a random closed set $\widetilde{\cal R} \subset [0,1]$.
Let $F_n$ be the size of the first part of ${\cal C}_n$, and
$L_n$ be the size of the last part of ${\cal C}_n$.
The following conditions are equivalent:

{\rm (i)} $\prob (F_n = 1) = \prob ( L_n = 1)$ for all $n$;

{\rm (ii)} $F_n \ed L_n$ for all $n$;

{\rm (iii)} $({\cal C}_n)$ is reversible; 

{\rm (iv)}  $({\cal C}_n)$ is symmetric;

{\rm (v)}  $({\cal C}_n)$ is the regenerative composition structure with
EPPF {\em \re{alpha-alpha}}, associated with an $(\alpha,\alpha)$ partition as 
in Section {\em \ref{alal}} for some $\alpha \in [0,1]$.
\end{theorem}

Before the proof of this result, we read from Theorem \ref{thm1} and
the discussion of Section \ref{alal} the following 
restatement of the equivalence of conditions (iii) and (v):

\begin{corollary}
For a random closed subset $\widetilde{\cal R}$ of $[0,1]$, the following
two conditions are equivalent:

{\rm (i)} $\widetilde{\cal R}$ is multiplicatively regenerative and 
$\widetilde{\cal R} \ed 1 - \widetilde{\cal R}$.

{\rm (ii)} $\widetilde{\cal R}$ is distributed like the zero set of a 
standard Bessel bridge of dimension $2 - 2 \alpha$, for some $\alpha \in [0,1]$.
\end{corollary}
\vskip0.5cm
\noindent
{\it Proof of Theorem {\rm \ref{symm}}.}
According to formula \re{qm}, for any regenerative composition structure
\eq
\label{firstone}
\prob( F_n = 1 ) = q(n:1)={ \OW ( n ) - \OW (n-1) \over \OW(n)/n }
\en
and the expressions \re{firstone} and \re{lastone} are obviously equal if $n = 1$ or $n = 2$.
We know that the $(\alpha,\alpha)$ regenerative composition structure is symmetric, hence reversible.
So for
$\OW_\alpha(n):= { [ 1 + \alpha ]_{n-1} / (n-1)! }$, the identity
$\prob( F_n = 1 ) = \prob(L_n = 1)$ together with \re{firstone} and \re{lastone}
yields 
\eq
\label{id1}
{ n \alpha \over n - 1 + \alpha} = n  - n {(-1)^n \over \OW_{\alpha}(n) } -  n \sum_{ k = 2 }^{n -1} { n - 1 \choose k - 1 } { (-1)^k \over \OW _\alpha (k) } .
\en
Suppose now that a regenerative composition structure is such that
$\prob( F_n = 1 ) = \prob(L_n = 1)$ for all $n = 1,2, \ldots$, and let us prove by
induction that its Laplace exponent $\OW$ normalised by $\OW(1) = 1$ is such that
\eq
\label{induc}
\OW( s ) = \OW_\alpha(s)
\en
for all $s = 1,2, \ldots$, where $\alpha \in [0,1]$ is defined by \re{induc} for $s = 2$, that is $\OW(2) = 1 + \alpha$.
According to \re{firstone} and \re{lastone}, we have for all $n = 2,3, \ldots$ that
\eq
\label{id2}
{ \OW ( n ) - \OW (n-1) \over \OW(n)/n } = n  - n {(-1)^n \over \OW(n) } -  n \sum_{ k = 2 }^{n -1} { n - 1 \choose k - 1 } { (-1)^k \over \OW (k) }
\en
so if we make the inductive hypothesis that \re{induc} holds for all $s \le n-1$ then we read from \re{id1} and \re{id2} that
$$
{ \OW ( n ) - \OW (n-1) \over \OW(n)/n } = { n \alpha \over n - 1 + \alpha} + n (-1)^n \left[ {1 \over \OW_\alpha(n) } - {1 \over \OW(n) } \right]
$$
which yields the expression
$$
\OW (n) = ( \OW_\alpha(n-1) - (-1)^n ) / ( 1 - \alpha(n-1 - \alpha) - (-1)^n/\OW_\alpha(n) ).
$$
But we know this formula holds for $\OW (n) = \OW_\alpha(n)$, so this must be the unique solution of the recursion,
and the inductive step is established.
Finally, the sequence $\OW(1), \OW(2), \ldots$ determines $\OW(s)$ for all $s \ge 0$, by consideration
of the second formula in \re{denom},
and the fact that a finite measure on $[0,1]$ is determined by its moments.
\endpf

\section{Transition probabilities}
\label{trans-probs}

Transition probabilities describing the succession of random compositions 
$({\cal C}_n)$ or ordered partitions $({\cal C}_n^*)$ as $n$ grows follow at once from the 
product formula \re{produ} for the composition probability function. 
For ordered partitions of $[n]$ these transition probabilities can be read immediately from
\re{pstar}, as indicated in James \cite[\S 5.4]{james03}.

\par Assuming that ${\cal C}_n^* = (A_1,\ldots,A_{\kell})$, 
an ordered partition  ${\cal C}_{n+1}^*$ of $[n+1]$ is 
obtained either by inserting singleton block $\{n+1\}$   into the sequence $A_1,\ldots,A_{\kell}$
or by adjoining the element $n+1$ to one of the blocks. 
It is easy to compute that  $n+1$ is inserted before $A_1$ 
with probability 
$$\frac{q(n+1:1)}{n+1}$$   
or adjoined to $A_1$  
with probability
$$\frac{\nlambda_{1}+1}{n+1}\,\frac{q(n+1:\nlambda_{1}+1)}{q(n:\nlambda_{1})}.$$
Inductively, with probability
$$
\prod_{i=1}^{j} 
\left( 1-\frac{q(\NLambda_i+1:1)}{\NLambda_i+1}-\frac{\nlambda_i+1}{\NLambda_i+1}\,\,\frac{q(\NLambda_i+1:\nlambda_i+1)}
{q(\NLambda_i:\nlambda_i)}
\right) 
$$
$n+1$ is neither inserted immediately before nor adjoined to one of the blocks $A_1,\ldots,A_{j}$,
and conditionally on this event (and given $(A_1,\ldots,A_{\kell})$) this element is  inserted as a singleton immediately 
following
$A_{j}$ with probability
$$\frac{q(\NLambda_{j+1}+1:1)}{\NLambda_{j+1}+1}$$
or  adjoined to $A_{j+1}$ (for $j<\kell$) with probability 
$$\frac{\nlambda_{j+1}+1}{\NLambda_{j+1}+1}\,\frac{q(\NLambda_{j+1}+1:\nlambda_{j+1}+1)}{q(\NLambda_{j+1}:\nlambda_{j+1})}.$$
Here, the $\nlambda_i$ are the sizes of the $A_i$ and the
$\NLambda_i$ are as in (\ref{produ}).

\par A transition law for integer compositions follows from the above. 
It is exactly the same as for the analogous ordered set partitions with the 
exception of the case when a composition of $n$
is changed by appending a 1 to a series of unit parts like $1,1,\ldots,1$,
in which case the transition probability is obtained by summation of individual probabilities 
of all possible singleton insertions into the series.

\section{Interval partitions}

The above probabilities of the two kinds  of transition (insertion and joining)
are equal to the expected sizes of intervals of a partition of $[0,1]$ induced by a uniform sample of $n$ points
and $\widetilde{\cal R}$. 
From this viewpoint, a better prediction of the `future'  compositions arising when more  points
are added  to the sample is obtained by conditioning on 
the actual sizes of intervals. 
\par At first  we shall describe a somewhat simpler distribution of the interval sizes for the $[0,\infty]$-partition, 
which can be seen as 
 discretisation of a subordinator in the spirit of \cite[Sections 3 and 4]{jp.bmpart}.
For each $n$,  a random set ${\cal R}$ and exponential order statistics 
$\epsilon_{1n},\ldots,\epsilon_{nn}$ induce a partition of 
$[0,\infty]$ associated with finite composition ${\cal C}_n$. 
The partition is comprised of
two kinds of parts: those containing some sample points or not.
The parts of the first kind are either
open interval components
of ${\cal R}^c$ which contain at least one of the $\epsilon_{jn}\,$'s, or one-point parts $\{\epsilon_{jn}\}$
corresponding to
$\epsilon_{jn}\in\cal R$ and appearing with positive probability only for $\drift>0$. The parts of the second kind
are the connected components (intervals or separate points) of the set resulting from removing  parts of the
first kind. The parts of different kinds interlace and if 
${\cal C}_n$ has $K_n$ classes there are $2K_n+1$ pieces of the partition, say 
$J_{1n},I_{1n},\ldots,J_{K_n-1,n},I_{K_n,n},J_{K_n+1,n}$,
which can be open or semiopen intervals or
one-point sets. 
Let $G_{1n},H_{1n},\ldots,G_{K_n-1,n},H_{K_n,n},G_{K_n+1,n}$ be the sizes of the parts, with slight
abuse of language we will call them `intervals', with understanding that some of them can degenerate into a point.

\vskip0.5cm

\begin{theorem} 
\label{t91}
The distribution of the random sequence $G_{1n},H_{1n},\ldots,G_{K_n-1,n},H_{K_n,n},G_{K_n+1,n}$ 
of interval sizes has the following properties:
\begin{itemize}
\item[{\rm (i)}] given the composition ${\cal C}_n$ all interval sizes are conditionally independent,
\item[{\rm (ii)}] $G_{1n}$ is independent of ${\cal C}_n$ and also independent of other interval sizes, and
has Laplace transform
\begin{equation}\label{G1}
{\ex }\exp\,(-s\,G_{1n})=\frac{\OW(n)}{\OW(n+s)}\,,
\end{equation}
\item[{\rm (iii)}] the unconditional distribution of $H_{1n}$ is given by 
\begin{equation}\label{H1}
{\prob}\,(H_{1n}\in {\rm d}z)=\frac{1-e^{-nz}}{\OW(n)}\, \Om({\rm d}z) + \frac{n\drift}{\OW(n)}\,\delta_0 ({\rm d}z),
\end{equation}
and given ${\cal C}_n$  the analogous
conditional distribution of $H_{1n}$ is 
$$\frac{{n\choose m}(1-e^{-z})^m\, e^{-(n-m)z} \,\Om({\rm d}z)+ n\,\drift\, 1(m=1)\,\delta_0({\rm d}z)}
{\ow(n:m)}$$
where $m$ is the first part of ${\cal C}_n$,

\item[{\rm (iv)}]  conditionally on the event that the first $j-1$ parts of ${\cal C}_n$   sum up to  $m$,
the truncated sequence $G_{jn},H_{jn},\ldots, H_{K_n,n},G_{K_n+1,n}$ is 
independent of  
the variables
$G_{1n}\,,H_{1n}\,,\ldots,G_{j-1,n}\,,H_{j-1,n}$ and of
the first $j-1$ parts of composition ${\cal C}_n$,  
and has the same distribution as the interlacing sequence 
$$G_{1,n-m}\,, H_{1,n-m}\,,\ldots,H_{K_{n-m}-j,n-m}\,,G_{K_{n-m}-j+1,n-m}$$
of interval sizes associated with the composition ${\cal C}_{n-m}$ of integer $n-m$.
\end{itemize}
\end{theorem}
{\it Proof.} The independence claims involved in (i) and (iv) follow from the memoryless property of the exponential distribution
and the strong Markov property of ${\widetilde{\cal R}}$ applied at the right endpoints of intervals $I_j$ or $J_j$.
Formulas (\ref{G1}), (\ref{H1}) follow  from Lemma \ref{firstpas} and the second formula in (iii) follows by routine 
conditioning.
\endpf

\par Mapping $[0,\infty]$ to $[0,1]$ by $z\mapsto 1-e^{-z}$ sends the partition of $[0,\infty]$ to
a partition of the unit interval, say 
$\widetilde{J}_{1n}\,, \widetilde{I}_{1n}\,,\ldots,\widetilde{I}_{K_n}\,,\widetilde{J}_{K_n+1}$,
which is the partition induced by a uniform sample and a multiplicatively regenerative set $\widetilde{\cal R}$.
The probability law of the partition of $[0,1]$ follows from Theorem \ref{t91}. Thus,
by virtue of the identity
${\ex }(1-\widetilde{G}_{1n})^s=\mathbb{E} \exp\,(-s\,G_{1n})$
the Laplace transform  (\ref{G1}) becomes a Mellin transform. 
Similarly,  the ratio $\widetilde{H}_{1n}/(1-\widetilde{G}_{1n})$ is independent of  $\widetilde{G}_{1n}$ and
has distribution
$${\prob}\left(\frac{\widetilde{H}_{1n}}{1-\widetilde{G}_{1n}}\in {\rm d}x\right)
=\frac{1-(1-x)^n}{\OW(n)}\,\,\widetilde{\Om}({\rm d}x)+ \frac{n\drift}{\OW(n)}\,\delta_0({\rm d}x)\,.$$
The distribution of the rest intervals follows recursively, by scaling with factor
$(1-\widetilde{G}_{1n}-\widetilde{H}_{1n})^{-1}$.

\par The sizes of these $2K_n+1$ intervals, say  $\widetilde{G}_{jn}$ and $\widetilde{H}_{jn}\,$,
determine the law of the extended composition when adding new sample points.
For example,
$${\ex }\widetilde{G}_{1n}=1-{\ex }\,(1-\widetilde{G}_{1n})=1-\frac{\OW(n)}{\OW(n+1)}=
\frac{\ow(n+1:1)}{(n+1)\OW(n+1)}$$
which by (\ref{qm})
is equal to 
$q(n+1:1)/(n+1)$ in accord with Section \ref{trans-probs}.
The sizes
also  
have a transparent  frequency interpretation  in terms of the
infinite composition ${\cal C}$.
For example, $\widetilde{G}_{1n}$ is the total 
frequency of the classes of ${\cal C}^*$ strictly preceding the first class represented in ${\cal C}^*_n$,
and $\widetilde{H}_{1n}$ is the frequency of the first class represented in ${\cal C}^*_n$.
\vskip0.5cm

\par {\bf Tripartite decomposition of $[0,1]$}
For $n=1$ the partition consists of three intervals $\widetilde{J}_{11}\,, \widetilde{I}_{11},
\widetilde{J}_{21}$ of sizes 
$G:=\widetilde{G}_{11}\,, H:=\widetilde{H}_{11},
D:=\widetilde{G}_{21}$.
The variable
$H$ is the frequency of the class of element $1$ and its distribution
is the structural distribution. Similarly, $G$ is the total frequency of
classes strictly preceding the class of 1 in ${\cal C}^*$, and $D$ is the total frequency of classes
strictly following the class of 1.
\par Moments of $G,H$ and $D$ 
have clear interpretation in terms of
finite compositions. Thus
\eq
\label{gone}
{\ex }\,(1-G)^{n-1}=\sum_{m=1}^n \frac{m}{n}\, q(n:m)=\frac{\OW(1)}{\OW(n)}
\en
is the probability that element 1 is in the first block of ${\cal C}_n^*$
or, what is the same, that a size-biased pick of a part from ${\cal C}_n$ yields the first part.
Similarly,
\eq
\label{gtwo}
{\ex }\,D^{n-1}=\frac{q(n:1)}{n}= \frac{\ow(n:1)}{n\,\OW(n)}
\en
is the probability that $\{1\}$  is the   first block  of ${\cal C}^*_n$.

\par  Furthermore,
the random variable
$H$ can be written as a product of two independent variables
$1-G$ and $H/(1-G)$, hence 
\eq
\label{gthree}
{\mathbb E}\left({H\over 1-G}\right)^{n-1}={{\mathbb E}H^{n-1}\over \mathbb E (1-G)^{n-1}}={\Phi(n:n)\over \Phi(1)}
\en
which is the conditional probability that the composition ${\cal C}_n^*$ is trivial given $1$ is in the first block.

\par For joint moments we have the formula
\begin{equation}\label{moments}
{\mathbb E}G^i\,H^{j-1}\,D^k=
\left(\sum_{a=0}^i {i\choose a}{(-1)^{a} \over \Phi(a+j+k)}\right)
\left(\sum_{b=0}^k (-1)^b{j\choose b}\Phi(j+b:j+b)\right)
\end{equation}
(the second sum may be further converted to variables $\Phi(1),\Phi(2),\ldots$) which follows 
from (\ref{gone}), (\ref{gthree}) and ${\mathbb E}H^n=p(n)=\Phi(n:n)/\Phi(n)$ by the binomial expansion of
 $$G^i\,H^j\,D^k= (1-(1-G))^j(1-G)^{j+k}\left({H\over 1-G}\right)^j \left(1-{H\over 1-G}\right)^k\,.$$
The joint moments have  the following interpretation.
Let $(A_1,A_2,A_3)$ be an ordered partition of $[n]$, $n=i+j+k$, such that $1\in A_2$ and 
the blocks are of sizes $i,j$ and $k$, respectively, with $i\geq 0$, $j\geq 1$ and $k\geq 0$.
Then (\ref{moments}) is the probability that $A_2$ is a block of ${\cal C}_n^*$ and 
$(A_1,A_2,A_3)$ is coarser than ${\cal C}_n^*$.
\par It follows that 
$${n-1\choose i,j-1,k} {\mathbb E}\, G^i\, H^{j-1}\, D^k$$ is  the probability that 
a size-biased pick of a part of ${\cal C}_n$ is $j$, and this part is preceded by a composition
of $i$ and followed by a composition of $k$ (with the obvious meaning when $i$ or $k$ is zero).
For $k=0$ this probability is equal to $(j/n){\mathbb P}(L_n=j)$ where $L_n$ is the last part of ${\cal C}_n$, 
computing this yields an  alternative proof for
(\ref{Ln}) and the formula for the Green matrix (\ref{Green}).

\paragraph{Acknowledgment}
Thanks to the referee for two careful readings of the paper, and for a number
of suggestions which helped to improve the exposition.

\def\cprime{$'$} \def\polhk#1{\setbox0=\hbox{#1}{\ooalign{\hidewidth
  \lower1.5ex\hbox{`}\hidewidth\crcr\unhbox0}}} \def\cprime{$'$}
  \def\cprime{$'$} \def\cprime{$'$}
  \def\polhk#1{\setbox0=\hbox{#1}{\ooalign{\hidewidth
  \lower1.5ex\hbox{`}\hidewidth\crcr\unhbox0}}} \def\cprime{$'$}
  \def\cprime{$'$} \def\polhk#1{\setbox0=\hbox{#1}{\ooalign{\hidewidth
  \lower1.5ex\hbox{`}\hidewidth\crcr\unhbox0}}} \def\cprime{$'$}
  \def\cprime{$'$} \def\cydot{\leavevmode\raise.4ex\hbox{.}} \def\cprime{$'$}
  \def\cprime{$'$} \def\cprime{$'$} \def\cprime{$'$}

\end {document}